\documentclass{amsart}
\usepackage{latexsym,amsxtra,amscd,ifthen}
\usepackage{amsfonts}
\usepackage{verbatim}
\usepackage{amsmath}
\usepackage{amsthm}
\usepackage{amssymb}

\numberwithin{equation}{section}

\theoremstyle{plain}
\newtheorem{theorem}{Theorem}[section]
\newtheorem{lemma}[theorem]{Lemma}
\newtheorem{proposition}[theorem]{Proposition}
\newtheorem{corollary}[theorem]{Corollary}

\theoremstyle{definition}

\newtheorem{examples}[theorem]{Examples}

\newtheorem{question}[theorem]{Question}

\makeatletter              
\let\c@equation\c@theorem  
\makeatother

\newcommand{\e}{{\textrm e}}

\newcommand{\rank}{\operatorname{rank}}

\newcommand{\boplus}{\textstyle \bigoplus}
\newcommand{\bcup}{\textstyle \bigcup}

\newcommand{\pp}{\mathfrak p}
\newcommand{\fg}{\mathfrak g}
\newcommand{\mm}{\mathfrak m}
\newcommand{\qq}{\mathfrak q}
\DeclareMathOperator{\Spec}{Spec}
\DeclareMathOperator{\height}{height}
\DeclareMathOperator{\gldim}{gldim}

\DeclareMathOperator{\Ext}{Ext}
\DeclareMathOperator{\Tor}{Tor}
\DeclareMathOperator{\pdim}{projdim}
\DeclareMathOperator{\fdim}{flatdim}

\DeclareMathOperator{\HB}{H}

\DeclareMathOperator{\ann}{ann}
\DeclareMathOperator{\lann}{l{.}ann}
\DeclareMathOperator{\rann}{r{.}ann}
\DeclareMathOperator{\cd}{cd}

\DeclareMathOperator{\gr}{gr}

\DeclareMathOperator{\Kdim}{Kdim}
\DeclareMathOperator{\rKdim}{rKdim}
\DeclareMathOperator{\Cdim}{Cdim}
\DeclareMathOperator{\injdim}{injdim}

\DeclareMathOperator{\GKdim}{GKdim}
\DeclareMathOperator{\End}{End}

\DeclareMathOperator{\D}{{\sf D}}
\DeclareMathOperator{\C}{{\sf C}}
\DeclareMathOperator{\Hom}{Hom}
\DeclareMathOperator{\RHom}{RHom}

\newcommand{\fm}{\mathfrak{m}}
\newcommand{\modleft}{\textup{-{\sf Mod}}}

\newcommand{\cal}{\mathcal}

\newcommand{\Oq}{{\mathcal O}_q}
\newcommand{\OqG}{\Oq(G)}
\newcommand{\CC}{{\mathbb C}}
\newcommand{\NN}{{\mathbb N}}
\newcommand{\QQ}{{\mathbb Q}}
\newcommand{\ZZ}{{\mathbb Z}}
\newcommand{\Ucheck}{{\check U}}
\newcommand{\tauhat}{\widehat{\tau}}
\newcommand{\deltahat}{\widehat{\delta}}
\newcommand{\tautil}{\widetilde{\tau}}
\newcommand{\deltil}{\widetilde{\delta}}
\newcommand{\qhat}{\widehat{q}}
\newcommand{\calC}{{\cal C}}
\DeclareMathOperator{\Fract}{Fract}
\DeclareMathOperator{\Htr}{Htr}
\newcommand{\Rtil}{\widetilde{R}}

\begin{document}

\title[Homological properties of $\OqG$]
{Homological properties of quantized coordinate rings of
semisimple groups}

\author{K.R. Goodearl and J.J. Zhang}

\address{Goodearl: Department of Mathematics, 
University of California at Santa Barbara,
Santa Barbara, CA 93106, USA}

\email{goodearl@math.ucsb.edu}

\address{Zhang: Department of Mathematics, Box 354350,
University of Washington, Seattle, WA 98195, USA}

\email{zhang@math.washington.edu}

\begin{abstract}
We prove that the generic quantized coordinate ring $\mathcal{O}_q(G)$
is  Aus\-land\-er-regular, Cohen-Macaulay, and catenary for 
every connected semisimple Lie group $G$. This answers 
questions raised by Brown, Lenagan, and the first author. We also prove
that under certain hypotheses concerning the existence of normal
elements, a noetherian Hopf algebra is Auslander-Gorenstein and
Cohen-Macaulay. This provides a new set of positive cases for a question
of Brown and the first author.
\end{abstract}

\subjclass[2000]{16A39, 16K40, 16E10, 16W50}





\keywords{Quantum group, quantized coordinate ring, Hopf algebra,
Auslander regular, Cohen-Macaulay, catenary}

\thanks{This research was partially supported by grants from the NSF
(USA) and by Leverhulme Research Interchange Grant F/00158/X (UK)}

\maketitle


\setcounter{section}{-1}
\section{Introduction}

A guiding principle in the study of quantized coordinate rings has been
that these algebras should enjoy noncommutative versions of the
algebraic properties of their classical analogs. Moreover, based on the
types of properties that have been established, one also conjectures
that quantized coordinate rings should enjoy properties similar to the
enveloping algebras of solvable Lie algebras. A property of the latter
type is the \emph{catenary} condition (namely, that all saturated chains
of prime ideals between any two fixed primes should have the same
length), which was established for a number of quantized coordinate rings
by Lenagan and the first author \cite{GL}. However, among the quantized
coordinate rings $\OqG$ for semisimple Lie groups $G$, catenarity has
remained an open question except for the case $G= SL_n$ \cite[Theorem
4.5]{GL}. The first goal of the present paper is to establish
catenarity for (the $\CC$-form of) all the algebras $\OqG$.

Following the principle indicated above, one expects the quantized
coordinate rings of Lie groups to be homologically nice; in the
noncommutative world, this means one looks for the
\emph{Auslander-regular} and \emph{Cohen-Macaulay} conditions. In fact,
Gabber's method for proving catenarity in enveloping algebras of
solvable Lie algebras relies crucially on these properties, as does
Lenagan and the first author's adaptation to quantum algebras. These
conditions were verified for $\Oq(SL_n)$ by Levasseur and Stafford, but
remained open for arbitrary $\OqG$, although Brown and the first author
were able to show that $\OqG$ has finite global dimension
\cite[Proposition 2.7]{BG1}. Our second goal here is to establish the
Auslander-regular and Cohen-Macaulay conditions for general $\OqG$.

Let $G$ be a connected semisimple Lie group over
$\mathbb C$. By $\OqG$, we mean the subalgebra of the Hopf dual of $U_q(\fg)$, where
$\fg$ is the Lie algebra of $G$, generated by the coordinate functions
of the type $\bf 1$ highest weight modules whose weights lie in the
lattice corresponding to a maximal torus of $G$. (See \cite[Definition
I.7.5]{BG2}, for instance.) We indicate that we consider the usual
one-parameter version of this algebra by calling it the \emph{standard}
quantized coordinate ring of $G$. Since some of our work relies on
properties established by Hodges-Levasseur-Toro
\cite{HLT} and Joseph
\cite{Jos}, we must restrict attention to the {\it complex form\/} of
$\OqG$, that is, we construct $U_q(\fg)$ and
$\OqG$ over the base field $\mathbb C$. The assumption that $q$ is {\it
generic\/} means that $q$ is a nonzero complex scalar which is not a
root of unity; in \cite{Jos}, $q$ is further restricted to be a scalar
transcendental over the rationals.

We summarize our main results for $\OqG$ as follows.

\begin{theorem}
\label{xxthm0.1} Let $G$ be a connected semisimple Lie group over
$\mathbb C$. Then the standard generic quantized coordinate ring
$\OqG$ \textup{(}defined over $\CC$\textup{)} is catenary, and Tauvel's
height formula holds. If $q$ is transcendental over $\QQ$, then $\OqG$ is
an Aus\-land\-er-regular and Cohen-Macaulay domain with $\gldim \OqG=
\GKdim
\OqG= \dim G$.  
\end{theorem}

The third goal of the paper concerns the question
\begin{enumerate}
\item[]
{\it If $H$ is a noetherian Hopf algebra, does $H$ have 
finite injective dimension?}
\end{enumerate}
This was asked by Brown and the first author \cite[1.15]{BG1} and
Brown \cite[3.1]{Br2}. Some partial results were given
by Wu and the second author when $H$ is PI \cite[Theorem 0.1]{WZ3}, and
when $H$ is graded with balanced dualizing complex 
\cite[Theorem 1]{WZ4}. Here we give another partial 
answer, for which we need to state some hypotheses. Let $A$ be a ring
(for the first condition) or an  algebra over a base field $k$ (for the
second and third). We say that $A$ satisfies
\begin{enumerate}
\item[(H1)]
 if $A$ is noetherian and $\Spec A$ is {\it normally 
separated\/}, namely, for every two primes $\pp\subsetneq \qq$ of $A$, 
there is a normal regular element in the ideal $\qq/\pp \subset A/\pp$.
\item[(H2)]
 if $A$ has an exhaustive ascending $\mathbb N$-filtration such that the
associated graded  ring $\gr A$ is connected graded noetherian with
{\it enough  normal elements\/}, namely, every non-simple graded prime
factor ring of
$\gr A$ contains a homogeneous normal element of positive degree.
\item[(H3)]
 if every maximal ideal of $A$ is co-finite 
dimensional over $k$.
\end{enumerate}
Our second result is the following.

\begin{theorem}
\label{xxthm0.2} Suppose $A$ is a Hopf algebra satisfying 
\textup{(H1,H2,H3)}. Then $A$ is Auslander-Gorenstein and 
Cohen-Macaulay and has a quasi-Frobenius classical quotient ring; 
further, $\Spec A$ is catenary and Tauvel's height formula holds.
\end{theorem}

It is known that the quantized coordinate rings $\OqG$ 
satisfy (H1) and (H2), but it is not known if they 
satisfy (H3). Therefore the proofs of Theorems 
\ref{xxthm0.1} and \ref{xxthm0.2} are similar in 
some steps and different in others.

Since our methods do not apply to quotient Hopf algebras of the $\OqG$,
we pose the following:

\begin{question}
\label{xxqn0.3} Let $G$ be a connected semisimple Lie group over
$\mathbb C$, and $I$ a Hopf ideal of the standard generic quantized
coordinate ring $\OqG$. Is the Hopf algebra $\OqG/I$
Auslander-Gorenstein and Cohen-Macaulay?
\end{question}

Throughout, let $k$ be a commutative base field. In treating general
properties of algebras, the base field will be $k$; in particular, the
unqualified term ``algebra'' will mean ``$k$-algebra'', and unmarked
tensor products will be $\otimes_k$. In working with $\OqG$, we will need
to assume $k$ is either $\mathbb C$ or a subfield of $\CC$.

Let $A$ be a ring or an algebra. Let $A^\circ$ denote the
opposite ring of $A$, and $A^\e$ (when $A$ is an algebra) the enveloping
algebra $A\otimes A^\circ$. Since we usually work with left modules, the
unqualified term ``module'' will mean ``left module''. Right
$A$-modules will be identified, when convenient, with
$A^\circ$-modules, and 
$(A,B)$-bimodules with $(A\otimes B^\circ)$-modules. Let $A\modleft$ and
$A^\circ\modleft$ denote the categories of all left and right
$A$-modules, respectively. 

When working with an $(A,B)$-bimodule $M$, the key property is that $M$
is finitely generated as both a left $A$-module and a right $B$-module
(note that this is stronger than requiring $M$ to be finitely generated
as a bimodule). We can record this property by saying that $M$ is {\it
left-right finitely generated\/}. Since all our rings $A$ and $B$ will be
noetherian, this is equivalent to saying that $M$ is {\it left-right
noetherian\/}, which we abbreviate to {\it noetherian\/}, following
common usage. To summarize: a {\it noetherian bimodule\/} is a bimodule
which is noetherian as both a left module and a right module.

\section{Consequences of (H1)}
\label{xxsec1}

Recall that a ring $A$ satisfies (H1) if $A$ is noetherian and 
$\Spec A$ is normally separated. It will be helpful to recall that these
hypotheses imply that $A$ satisfies the (left and right) strong second
layer condition
\cite[Theorem 12.17]{GW}. The main result of this section is  the
following proposition. Here, modules are treated as complexes
concentrated in degree 0, and we write $X\cong Y$ to indicate that
complexes $X$ and $Y$ are \emph{quasi-isomorphic}, that is, isomorphic
in the derived category.

\begin{proposition}
\label{xxprop1.1}
Suppose $A$ and $B$ are rings satisfying \textup{(H1)}. Let $X$ be a
bounded complex of $(A,B)$-bimodules such that the bimodule $\HB^i(X)$
is  noetherian for all $i$.
\begin{enumerate}
\item
If  $\Tor^A_i(A/\pp,X)=0$ for all $i\neq 0$ and all maximal 
ideals $\pp \subset A$, then $X\cong \HB^0(X)$ and $\HB^0(X)$ is a 
projective $A$-module.
\item
If, in addition, $\Tor^A_0(A/\pp,X)\neq 0$ for all maximal ideals 
$\pp \subset A$, then $\HB^0(X)$ is a progenerator in $A\modleft$. 
\end{enumerate}
\end{proposition}

We need several lemmas to prove this proposition.
Recall from \cite[Definition 1]{Br1} that a {\it bond\/} between two
noetherian  prime rings $A$ and $B$ is a nonzero noetherian
$(A,B)$-bimodule which is torsionfree both as a
left $A$-module and as a right 
$B$-module. The following lemma of the first author appears as
\cite[Theorem 10]{Br1}.

\begin{lemma}
\label{xxlem1.2} 
Let $A$ and $B$ be prime noetherian rings, and
$M$ a bond between them. If $B$ is simple, then $M$ is
projective over $A$. 
\end{lemma}

\begin{lemma}
\label{xxlem1.3} 
Let $A$ and $B$ be noetherian simple rings, and $M$ a
nonzero noetherian $(A,B)$-bimodule. Then $M$ is a progenerator  for
$A\modleft$ and for $B^\circ\modleft$.
\end{lemma}

\begin{proof} It follows from \cite[Lemma 8.3]{GW} that $M$ is
torsionfree on each side. Thus, $M$ is a bond, and so by Lemma
\ref{xxlem1.2},
$M$ is projective over 
$A$. Simplicity of $A$ then implies that $M$ is a generator over $A$.
Since  $M$ is finitely generated over $A$, it is a progenerator. By
symmetry,
$M$ is also a progenerator over $B$.
\end{proof}

Let $A$ be any ring. Let $N(A)$ denote the set of all regular normal
elements in $A$, and note that $N(A)$ is an Ore set. Then, let $S(A)$
denote the localization of 
$A$ by inverting all elements in $N(A)$. If the set of all
regular  elements in $A$ satisfies the left and right Ore conditions,
then we use $Q(A)$  to denote the localization of $A$ by inverting all
regular elements. The ring $Q(A)$ is called the {\it classical quotient
ring\/} of
$A$. Let $\Kdim$ denote the Gabriel-Rentschler Krull dimension. 

The following
lemma is clear.

\begin{lemma}
\label{xxlem1.4}
Suppose $A$ is a prime ring satisfying \textup{(H1)}.  Then every nonzero
ideal of $A$ contains a regular normal element.  Consequently, the 
localization $S(A)$ is simple.
\end{lemma}

\begin{lemma}
\label{xxlem1.5} 
Suppose $A$ and $B$ are prime rings satisfying
\textup{(H1)}, and that $M$ is a bond between them. 
\begin{enumerate}
\item
$S(A)\otimes_A M \cong M\otimes_B S(B) \cong S(A)\otimes_A M\otimes_B
S(B)$.
\item
$A$ is simple if and only if $B$ is simple.
\item
If $M\cong S(A)\otimes_A M$, then $A= S(A)$. In particular,
$A$ is simple.
\end{enumerate}
\end{lemma}

\begin{proof} (a) Since $M$ is a bond, left multiplication on $M$ by 
any $d\in N(A)$ is injective. Hence, $S(A)\otimes_A M$ is a directed
union of copies of $M$ as a right $B$-module. This implies that we can
make identifications
$$M\subset S(A)\otimes_A M\subset S(A)\otimes_A M\otimes_B S(B).$$

Let $e\in N(B)$, and let $N$ be the bimodule $M/Me$. Since $N$ is
torsion as a right $B$-module, it follows from \cite[Lemma 8.3]{GW} and
Lemma \ref{xxlem1.4} that there is an $f\in N(A)$ such that $f N=0$. 
This implies that
$S(A)\otimes_A N=0$, and hence $S(A)\otimes_A M=S(A)\otimes_A Me$. So,
the action of $e$ on $S(A)\otimes_A M$ is invertible.  This
says that
$S(A)\otimes_A M = S(A)\otimes_A M\otimes_B S(B)$. By symmetry,
$M\otimes_B S(B) = S(A)\otimes_A M\otimes_B S(B)$.

(b) Since $A$ and $B$ satisfy the second layer condition, they have the
same classical Krull dimension \cite[Corollary 14.5]{GW}. Part (b)
follows.

(c) Since $S(A)$ is simple and $M\cong S(A)\otimes_A M$ is a noetherian
$(S(A),B)$-bimodule, there is a positive integer $n$ such that $S(A)$
embeds in $M^{\oplus n}$ as a left $S(A)$-module \cite[Lemma 8.1]{GW}.
Consequently, $S(A)$ is a noetherian left $A$-module. This implies 
that, for each $e\in N(A)$, we have $Ae^{-i}=Ae^{-i+1}$ for some $i$.
Thus $Ae^{-1}=A$, and hence $e$ is invertible in $A$. Therefore
$S(A)=A$. Simplicity follows from Lemma \ref{xxlem1.4}.
\end{proof}

\begin{lemma}
\label{xxlem1.6} 
Suppose $A$ and $B$ are rings satisfying \textup{(H1)},
and $A$ is prime. Let $N$ be a nonzero noetherian $(A,B)$-bimodule.
\begin{enumerate} 
\item
If $S(A)\otimes_A N$ is nonzero, then it is a progenerator for
$S(A)\modleft$.  
\item
Similarly, if $A$ is simple, then $N$ is a
progenerator for 
$A\modleft$.
\item
If $N\cong S(A)\otimes_A N$, then $A=S(A)$ and $A$ is simple.
\end{enumerate}
\end{lemma}

\begin{proof} (a) First asssume that $S(A)\otimes_A N\ne  0$, and let
$T_1$ denote the $N(A)$-torsion sub(bi)module of $N$. Since
$S(A)\otimes_A N \cong S(A)\otimes_A (N/T_1)$, there is no loss of
generality in assuming that $N$ is $N(A)$-torsionfree. If $T_2$ is the
torsion sub(bi)module of $_AN$, then $\lann_A(T_2)$ is nonzero by
\cite[Lemma 8.3]{GW}. Then $\lann_A(T_2)$ contains an element of $N(A)$
by Lemma \ref{xxlem1.4}, and since $N$ is $N(A)$-torsionfree, we
conclude that $T_2=0$. Thus, $N$ is torsionfree as a left $A$-module.

Let $N_0=0 \subset N_1 \subset\cdots\subset N_m=N$
be a right affiliated series for $N$ (cf. \cite[Chapter 8]{GW}). The
ideals $\qq_i := \rann_B(N_i/N_{i-1})$ are prime, and $N_i/N_{i-1}$ is a
torsionfree right $(B/\qq_i)$-module as well as a torsionfree left
$A$-module \cite[Propositions 8.7, 8.9]{GW}. It suffices to show that
each $S(A)\otimes_A (N_i/N_{i-1})$ is a progenerator for
$S(A)\modleft$. Hence, we need only consider the case when $B$ is prime
and $N$ is a bond between $A$ and $B$.

Now we are in the situation of Lemma 
\ref{xxlem1.5}, that is, we have $N':=S(A)\otimes_A N\cong 
S(A)\otimes_A N\otimes_B S(B)$. Then $N'$ is a bond between
$S(A)$ and $S(B)$.  By Lemma \ref{xxlem1.4}, $S(A)$ and $S(B)$ are
simple; thus, the conclusion of part (a) follows from Lemma
\ref{xxlem1.3}.

(b) This is proved in the same manner as (a), with the help of
Lemma \ref{xxlem1.5}(b).

(c) The isomorphism $N\cong S(A)\otimes_A N$ implies that $N$ is
$N(A)$-torsionfree, and then, as in the proof of (a), it follows that
$_AN$ is torsionfree. Choose a right affiliated series for $N$, and let
$N'$ be its penultimate term. Then $\qq := \rann_B(N/N')$ is a prime
ideal of $B$, and $N/N'$ is a bond between $A$ and $B/\qq$. Since $N/N'$
is torsionfree as a left $A$-module, it is also $N(A)$-torsionfree, and
so $S(A)\otimes_A (N/N') \ne 0$. On the other hand, since the natural
map $N \rightarrow S(A)\otimes_A N$ is surjective, so is the natural map
$N/N' \rightarrow S(A)\otimes_A (N/N')$. Consequently, $N/N' \cong
S(A)\otimes_A (N/N')$, and therefore Lemma \ref{xxlem1.5}(c) implies that
$A=S(A)$.
\end{proof}

\begin{lemma}
\label{xxlem1.7}
Let $A$, $B$, $X$ satisfy the same hypotheses as \textup{Proposition
\ref{xxprop1.1}(a)}. Then $\Tor^A_i(N,X)=0$ for all noetherian
$A$-bimodules $N$ and all 
$i\neq 0$. As a consequence, $X\cong \HB^0(X)$. 
\end{lemma}

\begin{proof} Let $N$ be a noetherian $A$-bimodule. By induction on the
length of $X$, it is easy to see that each $\Tor^A_i(N,X)$ is a noetherian
$(A,B)$-bimodule.

We will use induction on $\rKdim N= \Kdim N_A$. If $r\Kdim N<0$, the
bimodule  $N=0$, and the assertion is trivial. Now assume that $\rKdim N
\ge0$, and that  $\Tor^A_i(N',X)=0$ for all $i\neq 0$ and all noetherian
$A$-bimodules  $N'$ with $\rKdim N'<\rKdim N$. We want to show that
$\Tor^A_i(N,X)=0$  for all $i\neq 0$.

By \cite[Corollary 8.8]{GW}, there exists a chain of sub-bimodules
$$N_0=0 \subset N_1 \subset\cdots\subset N_m=N$$
such that the annihilator ideals $\pp_j:= \lann_A(N_j/N_{j-1})$ and
$\qq_j :=\rann_A(N_j/N_{j-1})$ are prime and $N_j/N_{j-1}$ is a bond
between $A/\pp_j$ and $A/\qq_j$. It suffices to show that
$\Tor^A_i(N_j/N_{j-1},X) =0$ for all $i\neq 0$ and all $j$. Thus, there
is no loss of generality in assuming that $N$ 
is a bond between $A/\pp$ and $A/\qq$, for some prime ideals
$\pp,\qq\subset A$. If either $\pp$ or $\qq$ is a maximal ideal, then
both are maximal, by Lemma 
\ref{xxlem1.5}(b). Then, by Lemma \ref{xxlem1.2}, $N$ is 
projective over $A/\qq$. Since $\Tor^A_i(A/\qq,X)=0$ for 
all $i\neq 0$, it follows that $\Tor^A_i(N,X)=0$ for all $i\neq 0$. 

Now assume that $\pp$ and $\qq$ are both non-maximal, and let $A_1
:=A/\pp$ and $A_2 :=A/\qq$. For any $d\in N(A_1)$, one has $\rKdim
N/dN< \rKdim N$. Applying $\Tor^A_i(-,X)$ to the short exact sequence
$$0\to N \xrightarrow{l_d} N\to N/dN\to 0,$$
where $l_d$ is the left multiplication map by $d$, we obtain a long
exact sequence
$$\Tor^A_{i+1}(N/dN,X)\to \Tor^A_i(N,X)
\xrightarrow{l_d}\Tor^A_i(N,X)\to \Tor^A_i(N/dN,X)$$
for each $i$. 
For $i>0$ or $i<-1$, the two ends are zero by our induction hypothesis.
Hence, $T:=\Tor^A_i(N,X)$ is an $N(A_1)$-torsionfree 
$A_1$-module, and $T\cong S(A_1)\otimes_{A_1}T$. If $T\ne 0$, Lemma
\ref{xxlem1.6}(c) implies that $A_1$ is simple, a contradiction.
Therefore $T=0$ as  desired. 

For $i=-1$, we have a long exact sequence 
\begin{equation}
\Tor^A_{0}(N/dN,X)\to \Tor^A_{-1}(N,X)
\xrightarrow{l_d}\Tor^A_{-1}(N,X)\to \Tor^A_{-1}(N/dN,X).
\notag
\end{equation}
By the induction hypothesis, $\Tor^A_{-1}(N/dN,X)=0$, and so the map
$l_d$ is surjective. Let $U :=\Tor^A_{-1}(N,X)$, and let $V$ be the
$N(A_1)$-torsion submodule of $U$. Since $U$ is a bimodule,
so is $V$.  Now $l_d$ acts on $U/V$ bijectively. The argument for $T$
shows that
$U/V=0$, and so $U$ is $N(A_1)$-torsion. Since $U$ is a 
noetherian bimodule, there is a $d\in N(A_1)$ such that $d U=0$. But
$U=dU$, so $U=0$. Therefore $\Tor^A_i(N,X)=0$ for all $i\neq 0$, as
desired.

Finally, take $N=A$; then $\HB^i(X) \cong \Tor^A_i(A,X)=0$ for all $i\neq
0$.  Therefore $X\cong \HB^0(X)$.
\end{proof}

Now we are ready to prove Proposition \ref{xxprop1.1}.

\begin{proof}[Proof of Proposition \textup{\ref{xxprop1.1}}]
(a) By Lemma \ref{xxlem1.7}, $X$ is quasi-isomorphic to the bimodule
$\HB^0(X)$. Since $\Tor^A_i(A/\pp,\HB^0(X)) \cong \Tor^A_i(A/\pp,X) =0$
for all $i\ne 0$ and all maximal ideals $\pp\subset A$, we may replace
$X$ by $\HB^0(X)$. Thus, there is no loss of generality in assuming that
$X$ is a noetherian $(A,B)$-bimodule.

If $_AX$ is not projective, then it is not flat, and $\Tor^A_j(N,X)
\ne 0$ for some $j> 0$ and some finitely generated right $A$-module $N$.
By noetherian induction, we may assume that the ideal $\qq := \rann_A(N)$
is maximal for the existence of such an $N$, that is, $\Tor^A_i(M,X) =0$
for any $i\ne 0$ and any finitely generated right $A$-module $M$ whose
annihilator properly contains $\qq$. By a second noetherian induction, we
may assume that $N$ is a minimal criminal, that is, $\Tor^A_i(N/L,X) =0$
for any
$i\ne 0$ and any nonzero submodule $L\subseteq N$. If there is a nonzero
submodule $L\subseteq N$ with $\rann_A(L) \supsetneq \qq$, then
$\Tor^A_j(L,X) =0$ and $\Tor^A_j(N/L,X) =0$, by our two induction
hypotheses, which implies the contradiction $\Tor^A_j(N,X) =0$. Thus,
all nonzero submodules of $N$ have annihilator $\qq$, that is, $N$ is
fully faithful as a right $A/\qq$-module. It follows that $\qq$ is a
prime ideal of $A$ \cite[Proposition 3.12]{GW}. 

Set $A_2 := A/\qq$, and note that for any $c\in N(A_2)$, the set
$\ann_N(c)$ is a submodule of $N$, annililated by the nonzero ideal
$cA_2 \subset A_2$. This forces $\ann_N(c) =0$, since $N$ is a fully
faithful $A_2$-module. Thus, $N$ is $N(A_2)$-torsionfree. In particular,
the natural map 
$$N \longrightarrow N' := N\otimes_{A_2} S(A_2)= N\otimes_A S(A_2)$$
is injective. We identify $N$ with its image in $N'$.

Each finitely generated submodule $M\subset N'/N$ is annihilated by
some element of $N(A_2)$, whence $\rann_A(M) \supsetneq \qq$. By our
first noetherian induction, $\Tor^A_i(M,X) =0$ for all $i\ne 0$, for any
such $M$. Since $\Tor^A_i(-,X)$ commutes with direct limits, it follows
that $\Tor^A_i(N'/N,X) =0$ for all $i\ne 0$.

As a right $A$-module, $S(A_2)$ is the direct limit of the submodules
$c^{-1}A_2$ for $c\in N(A_2)$, and each $c^{-1}A_2 \cong A_2$. Since
$\Tor^A_i(A_2,X) =0$ for all $i\ne 0$, by Lemma \ref{xxlem1.7}, it
follows that $\Tor^A_i(S(A_2),X) = 0$ for all $i\ne 0$. Now set $K_{-1}
= X$, and choose short exact sequences of left $A$-modules
$$0 \to K_m \to F_m \to K_{m-1} \to 0$$
for all $m\ge 0$, where each $F_m$ is free. Since
$$\Tor^A_1(S(A_2),K_{m-1}) \cong \Tor^A_{m+1}(S(A_2),X) = 0$$
for all $m\ge 0$, the induced sequences
\begin{equation}
\label{E1.1}
0 \to S(A_2)\otimes_A K_m \to S(A_2)\otimes_A F_m \to S(A_2)\otimes_A
K_{m-1} \to 0
\tag{E1.1}
\end{equation}
are all exact. Further, $S(A_2) \otimes_A X \cong S(A_2)\otimes_{A/\qq}
(X/\qq X)$, which is a projective left
$S(A_2)$-module (possibly zero) by Lemma \ref{xxlem1.6}(b). Consequently,
we see that the short exact sequences \eqref{E1.1} all split. These
sequences remain split exact on ten\-sor\-ing with $N$, and hence we find
that the induced sequences
$$0 \to N'\otimes_A K_m \to N'\otimes_A F_m \to N'\otimes_A
K_{m-1} \to 0$$
are split exact. It follows that
$$\Tor^A_{m+1}(N',X) \cong \Tor^A_1(N',K_{m-1}) = 0$$
for all $m\ge 0$. 

From the long exact sequence for $\Tor$, we now have an exact sequence
$$0 = \Tor^A_{j+1}(N'/N,X) \to \Tor^A_j(N,X) \to \Tor^A_j(N',X) = 0.$$
However, this forces $\Tor^A_j(N,X) = 0$, contradicting our assumptions.
Therefore $X$ is indeed projective as a left $A$-module.

(b) As above, we may assume that $X$ is a noetherian bimodule. Let $I$ be
the trace ideal of $X$ in $A$, that is,
$$I= \sum_{f\in \Hom_A(X,A)} f(X).$$
For any maximal ideal $\pp\subset A$, we have $X/\pp X\ne 0$ by
hypothesis. Since $X$ is noetherian, $X/\pp X$ has a simple factor module,
and so there exists a maximal left ideal $\mm$ of $A$, containing
$\pp$, together with an epimorphism $f: X\rightarrow X/\pp X\rightarrow
A/\mm$ of left $A$-mod\-ules. Since $_AX$ is projective, $f$ lifts to an
$A$-module homomorphism $g: X\rightarrow A$ such that $g(X) \nsubseteq
\mm$. Consequently, $I \nsubseteq \pp$. 

Since $I$ is not contained in any maximal ideal of $A$, we conclude that
$I=A$. Therefore $X$ is a progenerator in $A\modleft$.
\end{proof}

We record the following analog of Proposition \ref{xxprop1.1}(a), in
which the hypothesis $\Tor^A_i(A/\pp,X)=0$ is only assumed for positive
indices $i$. Since we will not use the result in this paper, we omit the
proof, which is similar to that of Proposition \ref{xxprop1.1}(a). By
definition, the {\it flat dimension\/} of a complex
$X$ is
$$\fdim X := \sup\{i \mid \Tor^A_i(M,X)\neq 0 \text{ for all right
$A$-modules } M\}.$$

\begin{proposition}
\label{xxprop1.8}
Suppose $A$ and $B$ are rings satisfying \textup{(H1)}. Let $X$ be a
bounded complex of $(A,B)$-bimodules such that the bimodule $\HB^i(X)$
is  noetherian for all $i$. If  $\Tor^A_i(A/\pp,X)=0$ for all $i>0$ and 
 all maximal ideals $\pp \subset A$, then the flat dimension of 
$X$ is at most $0$. 
\end{proposition}

Finally, we note that $\OqG$ satisfies (H1).

\begin{theorem}
\label{xxlem1.9} Let $G$ be any connected semisimple Lie group
over ${\mathbb C}$. Then the standard generic quantized coordinate ring 
$\OqG$ is a noetherian domain of finite global dimension, and it
satisfies \textup{(H1)}.
\end{theorem}

\begin{proof} See \cite[Proposition 2.7]{BG1} and
\cite[Theorems I.8.9, I.8.18, and  II.9.19]{BG2}.
\end{proof}

\section{Consequences of (H2)}
\label{xxsec2} 

The main use of (H2) is to obtain the existence of an 
Auslander, Cohen-Macaulay, rigid dualizing complex. 
Since we will only use dualizing complexes as a tool in a 
few places, we intend to omit most of the definitions related
to dualizing complexes and derived categories. For those 
readers who would like to see more details about this topic,
we refer to the papers \cite{VdB, Ye1, YZ1, YZ2}. However,
we will review below the definitions of the Auslander 
and Cohen-Macaulay properties, since these two properties are 
the main point of Theorems \ref{xxthm0.1} and \ref{xxthm0.2}.

Recall that a {\it connected graded\/} algebra is an $\mathbb N$-graded
$k$-algebra $A= \bigoplus_{i\ge 0} A_i$ such that $A_0=k$ and $\dim_k
A_i <\infty$ for all $i$. A noetherian connected graded algebra
$A$ is said  to {\it have enough normal elements\/} if every non-simple
graded prime factor $A/\pp$ has a homogeneous normal element of positive 
degree. (Such a normal element is necessarily regular.) Recall that an
algebra
$A$ satisfies (H2) if there is an ascending $\mathbb N$-filtration
$F=\{F_iA\}_{i\geq 0}$ by $k$-subspaces such that
\begin{enumerate}
\item[(i)] $1\in F_0A$, 
\item[(ii)] $A=\bcup_{i\geq 0} F_i A$, 
\item[(iii)] $(F_iA)(F_jA)\subset F_{i+j}A$ for all $i,j$, and 
\item[(iv)] the 
associated graded ring $\gr A:=\boplus_{i=0}^{\infty} 
F_iA/F_{i-1}A$ is connected graded noetherian and has 
enough normal elements.
\end{enumerate}
Condition (iv) implies that $\gr A$ is a finitely generated algebra with
finite GK-dim\-en\-sion \cite[Proposition 0.9]{YZ1}, and consequently $A$
satisfies those same properties.

In this section, we will recall several concepts which 
are essential to our work. First, let $A$ be a noetherian ring. We say
that $A$ is {\it Gorenstein\/} if it has finite injective
dimension on both sides. The standard Auslander-Gorenstein,
Auslander-regular, and Cohen-Macaulay conditions are usually defined in
terms of grades of modules, as follows. For any finitely generated left
or right
$A$-module, the {\it grade\/} or the {\it $j$-number\/} of $M$ {\it
with respect to $A$\/} is defined to be
$$j(M) := \inf\{n\;|\; \Ext^n_A(M,A)\neq 0\}.$$
The ring $A$ is called {\it Auslander-Gorenstein\/} if it is Gorenstein
and it satisfies the {\it Auslander condition\/}:
\begin{enumerate}
\item[] For every finitely generated left (respectively, right)
$A$-module $M$ and every positive integer $q$, one has $j(N) \ge q$ for
every finitely generated right (respectively, left) $A$-sub\-mod\-ule
$N\subseteq
\Ext^q_A(M,A)$.
\end{enumerate}
An {\it Auslander-regular\/} ring is a noetherian, Auslander-Gorenstein
ring which has finite global dimension. The final condition requires $A$
to be an algebra with finite GK-dimension. We say that $A$ is {\it
Cohen-Macaulay\/} ({\it with respect to $\GKdim$\/}) if 
$$j(M)+\GKdim M= \GKdim A$$
for every nonzero finitely generated left or right $A$-module $M$.

These conditions have analogs for rigid dualizing complexes. Thus, let
$R$ be a rigid dualizing complex over $A$, and let $M$ be a finitely
generated left or right $A$-module. The {\it grade\/} of $M$ {\it with
respect to $R$\/} is defined to be
$$j_R(M) := \inf\{n\;|\; \Ext^n_A(M,R)\neq 0\}.$$
We say that $R$ is {\it Auslander\/} if
\begin{enumerate}
\item[]
For every finitely generated left (respectively, right) $A$-module $M$,
every integer $q$, and every finitely generated right (respectively, left)
$A$-submodule $N\subseteq \Ext^q_A(M,R)$, one has $j_R(N)\geq q$.
\end{enumerate}
Now suppose that $A$ is an algebra with finite GK-dimension.
We say that $R$ is {\it Cohen-Macaulay\/} ({\it with
respect to $\GKdim$\/}) if 
$$j_R(M)+\GKdim M=0$$
for every nonzero finitely generated left or right $A$-module $M$.
 
If $A$ is Gorenstein, then any shift $A[n]$ is a dualizing complex over
$A$. This leads to the following connection between the two versions of
the above concepts, which is easily checked once one notes that $j_R(M)=
j(M)-n$ for all finitely generated left or right $A$-modules $M$.

\noindent {\it A noetherian
algebra
$A$ with
$\GKdim A = n <\infty$ is Auslander-Gorenstein and Cohen-Macaulay if and
only if\/}
\begin{enumerate}
\item[(i)] {\it $A$ is Gorenstein, and\/}
\item[(ii)] {\it the dualizing complex $A[n]$ is Auslander and
Cohen-Macaulay.\/}
\end{enumerate}

\begin{proposition} 
\label{xxlem2.1} 
Suppose $A$ is an algebra satisfying \textup{(H2)}. 
\begin{enumerate}
\item $A$ has an 
Auslander, Cohen-Macaulay, rigid dualizing complex $R$. 
As a consequence, the Gelfand-Kirillov dimension of any finitely
generated $A$-module is an integer. 
\item For any ideal $I\subset A$, the complex $\RHom_A(A/I,R)$ is an 
Auslander, Cohen-Macaulay, rigid dualizing complex over $A/I$.
\end{enumerate}
\end{proposition}

\begin{proof} Part (a) is \cite[Corollary 6.9(iii)]{YZ1}, while part (b)
follows from \cite[Theorem 6.17, Proposition 3.9]{YZ1}.
\end{proof}

Given an algebra $A$ with a dualizing complex $R$, define functors  $D$
and $D^\circ$ so that
$$D(M)=\RHom_A(M,R)\qquad \text{and}\qquad
D^\circ(N)=\RHom_{A^\circ}(N,R)$$
for $M\in \D(A\modleft)$ and $N\in \D(A^\circ\modleft)$, where $\D(\C)$
denotes the derived category of complexes over an abelian category $\C$.
These functors provide dualities between the subcategories
$\D_f(A\modleft)$ and $\D_f(A^\circ\modleft)$ whose objects are the
complexes with finitely generated cohomology modules \cite[Proposition
1.3]{YZ1}. Let
$(-)'=
\Hom_k(-,k)$ denote the $k$-linear vector space dual.

We shall need the following relatively standard fact: \emph{Suppose
that $A$ and $B$ are noetherian rings, $X$ a bounded complex of left
\textup{(}respectively, right\textup{)}
$A$-modules, and $Y$ a bounded complex of $(A,B)$-bimodules
\textup{(}respectively, $(B,A)$-bimodules\textup{)}. If the cohomology
modules $\HB^i(X)$ are finitely generated and the bimodules
$\HB^i(Y)$ are noetherian for all
$i$, then the right \textup{(}respectively, left\textup{)}
$B$-modules
$\Ext^i_A(X,Y)$ \textup{(}respectively,
$\Ext^i_{A^\circ}(X,Y)$\textup{)} are finitely generated.} By induction on
the lengths of the complexes
$\HB^*(X)$ and $\HB^*(Y)$, this can be reduced to the case where $X$ and
$Y$ are (bi)modules. That case is well known.

\begin{lemma}
\label{xxlem2.2} Suppose $A$ is an algebra satisfying \textup{(H2)}, and 
let $R$ be the dualizing complex given in 
\textup{Proposition \ref{xxlem2.1}(a)}. Set $D(-) := \RHom_A(-,R)$ as
above.
\begin{enumerate}
\item $D(M)\cong M'$ for all finite dimensional $A$-modules
$M$.
\item $\Ext^i_A(A/I,A)$ is a noetherian $(A/I,A)$-bimodule for all ideals
$I\subset A$ and all $i\ge 0$.
\end{enumerate}
\end{lemma}

\begin{proof} (a) Let $M$ be a finite dimensional $A$-module and $\pp
:=\ann_A(M)$. Then
$A/\pp$ is finite dimensional over $k$. One has
$$D(M)\cong \RHom_A((A/\pp)\otimes^L_{A/\pp} M,R)
\cong \RHom_{A/\pp}(M, D(A/\pp)) =: (*)$$
\cite[Theorem 10.8.7]{Wbl}.
By Proposition \ref{xxlem2.1}(b) (or \cite[Theorem 3.2]{YZ1} in this
special case),
$R_{A/\pp}:=D(A/\pp)$ is a rigid  dualizing complex over $A/\pp$. Since
$A/\pp$ is finite dimensional over
$k$, we have
$R_{A/\pp}\cong (A/\pp)'$. Hence
$$(*)\cong \RHom_{A/\pp}(M, (A/\pp)')\cong M'.$$

(b) Let  $I$ be an ideal of $A$ and $i$ a nonnegative integer. By the
comments above, $\Ext^i_A(A/I,A)$ is finitely generated as a right
$A$-module. By Proposition \ref{xxlem2.1}(b), $D(A/I)$ is a rigid
dualizing complex over $A/I$. By duality,
$$\Ext^i_A(A/I,A) \cong \Ext^i_{A^\circ}(D(A),D(A/I)) \cong
\Ext^i_{A^\circ}(R,D(A/I)).$$
The latter bimodule is finitely generated as a left $(A/I)$-module by the
comments above.
\end{proof}

Part (a) of the lemma says that if $M$ is an $A$-module finite
dimensional  over $k$, then 
\begin{equation}
\Ext^i_A(M,R)=\begin{cases} M'
& {\text {if}} \quad i=0,\\
0 & {\text {if}} \quad i\neq 0.
\end{cases}
\notag
\end{equation}

Let $\D^b_{f}(A\modleft)$ denote the derived category of bounded
complexes $Y$ of $A$-mod\-ules such that $\HB^i(Y)$ is finitely generated 
over $A$ for all $i$.

\begin{lemma}
\label{xxlem2.3}
Suppose $A$ is an algebra satisfying \textup{(H2)}. Let $X$ and $Y$ be
objects in 
$\D^b_f(A\modleft)$ such that $\HB^i(X)$ is finite dimensional
over $k$ for all $i$. Then $\Ext^i(X,Y)$ is finite dimensional
over $k$ for all $i$.
\end{lemma}

\begin{proof} By induction on the length of $X$, we may assume
that $X$ is a finite dimensional $A$-module $M$. By duality,
$\Ext^i_A(M,Y)\cong \Ext^i_{A^\circ}(D(Y),D(M))$. By Lemma
\ref{xxlem2.2}(a), $D(M)\cong M'$, which is finite dimensional;
by the properties of the duality $D$, we have $D(Y)\in 
\D^b_f(A^\circ\modleft)$. So, it remains to show that
$\Ext^i_{A^\circ}(Z,N)$ is  finite dimensional when $Z\in
\D^b_f(A^\circ\modleft)$ and
$N$ is a finite dimensional right $A$-module. But this
follows by taking a bounded above projective resolution $P$
of $Z$ such that each term of $P$ is finitely generated.
\end{proof}

A ring $A$, or its prime spectrum $\Spec A$, is called {\it catenary\/}
if, for every pair of prime ideals $\pp\subset \qq$, all saturated
chains of primes between
$\pp$ and $\qq$ have same length. We say that {\it Tauvel's 
height formula\/} holds in an algebra $A$ provided
$$\height \pp +\GKdim A/\pp=\GKdim A$$
for all prime ideals $\pp \subset A$.
These properties are consequences of (H1) and (H2), as follows.

\begin{theorem}
\cite[Theorems 0.1 and 2.23]{YZ1}
\label{xxlem2.4}
Suppose $A$ is an algebra satisfying \textup{(H1,H2)}.
\begin{enumerate}
\item
$\Spec A$ is catenary. 
\item
If, in addition, $A$ is prime, then Tauvel's height
formula holds.
\end{enumerate}
\end{theorem}

Even if $A$ is not prime, but rather equidimensional in a suitable
sense, Tauvel's height formula still holds.

\begin{corollary}
\label{xxlem2.5}
Suppose $A$  is an algebra satisfying \textup{(H1,H2)}. If $\GKdim A/\pp
=\GKdim A$ for all minimal primes $\pp$, then Tauvel's height formula
holds.
\end{corollary}

\begin{proof} Let $\qq$ be any prime, and let $\pp$ be a minimal
prime such that $\pp \subset \qq$ and the length of a saturated chain
between $\pp$ and
$\qq$ is equal to the height of $\qq$. Since $A/\pp$ satisfies 
\textup{(H1)} and \textup{(H2)} also, Theorem \ref{xxlem2.4}(b)
implies that
$$\height \qq/\pp+\GKdim A/\qq=\GKdim A/\pp.$$
By the choice of $\pp$, one has $\height \qq/\pp=\height \qq$. By 
hypothesis, $\GKdim A/\pp=\GKdim A$. Therefore, 
$$\height \qq +\GKdim A/\qq=\GKdim A.$$
\end{proof}

Next, we consider the quantized coordinate rings $\OqG$.

\begin{theorem}
\label{xxprop2.6} 
If $G$ is any connected semisimple Lie group
over ${\mathbb C}$, then the standard generic quantized coordinate ring
$\OqG$ satisfies 
\textup{(H2)}. 
As a consequence, $\Spec \OqG$ is 
catenary and Tauvel's height formula holds.
\end{theorem}

\begin{proof} By the proof of \cite[Theorem I.8.18]{BG2},
$A:=\OqG$ has a set of generators satisfying the hypotheses of 
\cite[Proposition I.8.17]{BG2}. By the proof of 
\cite[Proposition I.8.17]{BG2}, $A$ has an $\mathbb N$-filtration 
such that $\gr A$ is a graded factor ring of a 
skew polynomial ring 
$${\cal O}_{\mathbf q}({\mathbb C}^m) := {\mathbb C}\langle x_1,\dots,x_m
\mid x_ix_j= q_{ij}x_jx_i \text{ for all } i,j \rangle,$$
for some matrix ${\mathbf q}= (q_{ij})$ of nonzero scalars,
with $\deg x_i>0$ for each $i$. It is clear that
${\cal O}_{\mathbf q}({\mathbb C}^m)$ is connected graded noetherian, and
that it has enough normal elements (namely, the images of $x_1,\dots,x_m$
in graded factor rings). Any graded factor ring of
${\cal O}_{\mathbf q}({\mathbb C}^m)$ inherits these properties. Hence,
$A$ satisfies (H2).

The second assertion follows from Theorems
\ref{xxlem1.9} and \ref{xxlem2.4}.
\end{proof}

Here is one final 
consequence of (H2). Recall that an algebra $A$ is called {\it
universally noetherian\/} if $A\otimes B$ is
noetherian for every noetherian algebra $B$ \cite[p.596]{ArSZ}.

\begin{proposition}
\label{xxlem2.7}
If $A$ is an algebra satisfying \textup{(H2)}, then it is universally 
noetherian.
\end{proposition}

\begin{proof} This follows from
\cite[Propositions 4.3, 4.9 and 4.10]{ArSZ}.
\end{proof}

\section{relatively projective modules}
\label{xxsec3}

The hypothesis (H3) turns out to be very useful in dealing with modules
satisfying the properties described in the following lemma. We leave the
proof to the reader.

\begin{lemma}
\label{xxlem3.4}
Let $A$ be an algebra and $M$ a finitely generated $A$-module. Then the
following are equivalent:
\begin{enumerate}
\item
$M/\pp M$ is projective over $A/\pp$ for all co-finite 
dimensional ideals $\pp \subset A$.
\item
$\Hom_A(M,-)$ is an exact functor on finite dimensional
$A$-modules.
\item
$-\otimes_A M$ is an exact functor on finite dimensional right
$A$-modules.
\end{enumerate}
\end{lemma}

We shall say that a finitely generated module $M$ over an algebra $A$ is
{\it relatively projective\/} (with respect to the category of finite
dimensional $A$-modules) provided $M$ satisfies the equivalent
conditions of Lemma \ref{xxlem3.4}. 

The hypothesis (H3), which says that every maximal ideal is co-finite
dimensional, seems quite restrictive. Nonetheless, it holds in a number
of interesting non-PI situations. Here are some examples.

\begin{examples}
\label{xxex3.1} {\bf (a)} Let $A$ be a finitely generated prime algebra.
If $A$ is not simple and $\GKdim A\leq 2$, then $A$ satisfies
\textup{(H3)}. To see this, we let $\pp$ be a maximal ideal of $A$. Since
$A$ is not simple, $\pp\neq 0$, and hence $\GKdim A/\pp\leq 1$. By
\cite{SSW}, $A/\pp$ is PI. Now simple affine PI rings are finite
dimensional, so $\pp$ is co-finite dimensional.

{\bf (b)} If $A$ is a finitely generated prime Hopf algebra with
$\GKdim A\le 2$, then $A$  satisfies (H3). This follows from (a), since
$A$ cannot be simple unless the kernel of its counit is zero, in which
case $\dim_k A = 1$. Of course, not all Hopf algebras satisfy (H3), as
witnessed by many group algebras and enveloping algebras.

{\bf (c)} If $k$ is algebraically closed and $q\in k^\times$ is not a
root of unity, then the quantized coordinate rings $\Oq(k^2)$ and
$\Oq(SL_2(k))$ satisfy (H3) -- just refer to the pictures of their prime
spectra in \cite[Diagrams II.1.2, II.1.3]{BG2}.

{\bf (d)} Assume again that $k$ is algebraically closed and $q\in
k^\times$ is not a root of unity. We generalize the first example of (c)
to the algebra
$$A= \Oq(k^n) := k\langle x_1,\dots,x_n \mid x_ix_j= qx_jx_i \text{\ for
all\ } i<j \rangle.$$
Given a maximal ideal $\pp$ in $A$, let $I= \langle x_i \mid x_i\in
\pp\rangle$. It suffices to show that the maximal ideal $\pp/I$ in $A/I$
is co-finite dimensional, and $A/I\cong \Oq(k^m)$ where $m= n-
\text{card}(I)$. Hence, there is no loss of generality in assuming that
$x_i \notin \pp$ for all $i$.

Now the images of the $x_i$ in $A/\pp$ are nonzero normal elements,
hence invertible by simplicity. Thus, $A/\pp \cong B/\pp B$ where
$$B= \Oq((k^\times)^n) := A[x_1^{-1},\dots,x_n^{-1}].$$
It is known that the maximal ideals of $B$ are induced from the maximal
ideals of its center \cite[Corollary 1.5]{GLet2}. This center can be
calculated using \cite[Lemma 1.2]{GLet2} -- if $n$ is even, then
$Z(B)=k$, while if $n$ is odd, then $Z(B)= k[z^{\pm1}]$ where $z :=
x_1x_2^{-1}x_3x_4^{-1} \cdots x_{n-1}^{-1}x_n$.

In case $n$ is even, we obtain $\pp B=0$ and $\pp=0$, meaning that $A$
itself is simple. This is only possible when $n=0$, that is, $A=k$, and
thus $\dim_k A/\pp =1$ in this case.

In case $n$ is odd, $\pp B= (z-\alpha)B$ for some $\alpha\in k^\times$,
and thus the element
$$w := x_1x_3\cdots x_n - \alpha x_2x_4\cdots x_{n-1}$$
lies in $\pp B\cap A= \pp$. Since $w$ generates $\pp B$, we have
$$\pp= \{a\in A\mid (x_1x_2\cdots x_n)^ta \in wA \text{\ for some\ }
t\in \NN\},$$
from which it is easy to check that, in fact, $\pp= \langle w\rangle$.
Now if $n\ge 3$, it would follow that $\pp \subseteq \langle
x_1,x_2\rangle$, contradicting the maximality of $\pp$. Thus, $n=1$,
whence $\pp= \langle x_1-\alpha\rangle \subset A= k[x_1]$, and again
$\dim_k A/\pp =1$.

Therefore all maximal ideals of $A$ have codimension 1.

{\bf (e)} In contrast to (d), multiparameter quantum affine spaces need
not satisfy (H3). For instance, let $q\in k^\times$ be a non-root of
unity, and let
$$A := k\langle x,y,z \mid xy=qyx,\, xz= q^{-1}zx,\, yz= qzy \rangle.$$
Then $xyz\in Z(A)$, and $A/ \langle xyz-1 \rangle \cong
\Oq((k^\times)^2)$, which is a simple algebra \cite[Example
1.8.7(ii)]{MR}. Thus, $\langle xyz-1 \rangle$ is a maximal ideal of $A$
with infinite codimension.  \hfill$\square$
\end{examples}

\begin{question}\label{xxqn3.2}
Let $G$ be any connected semisimple Lie group over $\CC$. Do all maximal
ideals of the standard generic quantized coordinate ring $\OqG$ have
finite codimension?
\end{question}

 The main result in this 
section is the following proposition.

\begin{proposition}
\label{xxprop3.2} Suppose that $A$ is an algebra satisfying
\textup{(H1,H3)}, $B$ a ring satisfying \textup{(H1)}, and
$M$ a noetherian $(A,B)$-bimodule. If $_AM$ is relatively projective, 
then $_AM$ is projective. If, in addition, $M/\pp M\neq 0$
\textup{(}or, equivalently, $\Hom_A(M,A/\pp)\neq 0$\textup{)} for 
all maximal ideals $\pp$ of $A$, then $_AM$ is a progenerator for
$A\modleft$.
\end{proposition}

We need some lemmas. Recall that a module $M$ over an algebra $A$ is
said to be {\it locally finite dimensional\/} if every finitely generated
submodule of $M$ is finite dimensional.

\begin{lemma}
\label{xxlem3.3} Suppose $A$ is an algebra satisfying \textup{(H1)}. 
Then the injective hull of any finite dimensional left or right
$A$-module is locally finite dimensional.
\end{lemma}

\begin{proof}
Since the lemma is left-right symmetric, it suffices to prove it for
right modules.

Let $E$ be the injective hull of a finite dimensional right $A$-module
$M$. Then $E= E_1\oplus \cdots\oplus E_m$ for some uniform injective
modules $E_j$, which are injective hulls of the finite dimensional
modules $E_j\cap M$. Since it suffices to show that each $E_j$ is
locally finite dimensional, there is no loss of generality in assuming
that $E$ is uniform. Now $M$ has an essential simple submodule $M'$, and
we may replace $M$ by $M'$, so we may assume that $M$ is simple. As a
consequence, $\pp := \rann_A(M)$ is a co-finite dimensional maximal
ideal of $A$. Moreover, $\pp$ is the \emph{assassinator} of $E$, that
is, the unique associated prime.

Recall from Section \ref{xxsec1} that (H1) implies the strong second
layer condition. By \cite[Corollary 12.8]{GW}, each finitely generated
submodule of
$E$ is annihilated by a product of primes from the right link closure of
$\pp$. All such primes must have finite codimension, because finite
codimension carries across links: if $\pp_1$ and $\pp_2$ are prime
ideals of $A$ such that $\pp_1 \rightsquigarrow \pp_2$ via
$(\pp_1\cap\pp_2)/I$, then $A/\pp_1$ is finite dimensional if and only
if $(\pp_1\cap\pp_2)/I$ is finite dimensional, if and only if $A/\pp_2$
is finite dimensional. Thus, the annihilator of each finitely generated
submodule of
$E$ has finite codimension. Therefore $E$ is locally finite dimensional.
\end{proof}

\begin{lemma}
\label{xxlem3.5} Suppose $A$ is an algebra satisfying \textup{(H1)} 
and $M$ a finitely generated relatively projective $A$-module.
\begin{enumerate}
\item
$\Ext^i_A(M,N)=0$ for all $i>0$ and all finite dimensional
$A$-modules $N$.
\item
$\Tor_i^A(S,M)=0$ for all $i>0$ and all finite dimensional right
$A$-modules $S$.
\end{enumerate}
\end{lemma}

\begin{proof} (a) Suppose 
$$0\to N\to I^0\to I^1\to I^2\to \cdots$$ 
is a minimal injective resolution of $N$. The assertion follows if the
complex
$$\Hom_A(M,I^0)\to \Hom_A(M,I^1)\to 
\Hom_A(M,I^2)\to \cdots$$
is exact. Since each $I^i$ is locally finite
dimensional by Lemma \ref{xxlem3.3}, it suffices to show 
that $\Hom_A(M,-)$ is exact on locally finite dimensional modules. 
Suppose
$$0\to P \xrightarrow{f} Q \xrightarrow{g} R\to 0$$
is a short exact sequence of locally finite dimensional $A$-modules,
and let $h:M\to R$ be a homomorphism. Then $R'= h(M)$ is a finitely
generated submodule of
$R$, and we can choose a finitely generated submodule $Q' \subseteq Q$
such that
$g(Q')= R'$. Since $Q'$ is finite dimensional, we obtain a short exact
sequence
$$0\to P' \xrightarrow{f'} Q' \xrightarrow{g'} R'\to 0$$
of finite dimensional $A$-modules, where $P' := f^{-1}(Q')$
and $f'$, $g'$ are the appropriate restrictions of $f$, $g$. When $h$ is
viewed as a homomorphism $M\to R'$, it lifts to a homomorphism $h': M\to
Q'$ because $M$ is relatively projective. This provides the required
lifting of $h$ to $\Hom_A(M,Q)$, and proves that the sequence
$$0\to \Hom_A(M,P)  \xrightarrow{f^*} \Hom_A(M,Q)  \xrightarrow{g^*}
\Hom_A(M,R) \to 0$$ 
is exact, as desired.

(b) Let $S$ be an arbitrary finite dimensional right $A$-module. Since
$k$ is injective as a module over itself, \cite[Proposition VI.5.1]{CE}
implies that
$$\Hom_k \bigl( \Tor^A_i(S,M),k \bigr) \cong \Ext^i_A \bigl(
M,\Hom_k(S,k) \bigr)$$ for all $i$. By part (a), $\Ext^i_A \bigl(
M,\Hom_k(S,k) \bigr) =0$ for $i>0$, and therefore we obtain
$\Tor^A_i(S,M) =0$ for
$i>0$.
\end{proof}

Now we are ready to prove Proposition \ref{xxprop3.2}.

\begin{proof}[Proof of Proposition \textup{\ref{xxprop3.2}}]
Let $M$ be a noetherian $(A,B)$-bimodule which is relatively projective
as an $A$-module.  Since every maximal ideal of $A$ is co-finite
dimensional by (H3), Lemma \ref{xxlem3.5}(b) implies that
$\Tor_i^A(A/\pp,M)=0$ for all
$i\ne 0$ and all maximal ideals $\pp$ of $A$. By Proposition
\ref{xxprop1.1}(a),
$_AM$ is projective. If, in addition, $M/\pp M\neq 0$ for all maximal
ideals $\pp$ of $A$, then $_AM$ is a progenerator by Proposition
\ref{xxprop1.1}(b).
\end{proof}

\section{Proof of Theorem \ref{xxthm0.2}}
\label{xxsec4}

An algebra $A$ is said to satisfy
\begin{enumerate}
\item[(H4)] if for all nonnegative integers $n$,
\begin{enumerate}
\item the functor $\Ext^n_A(-,A)$ is exact on finite
dimensional $A$-modules;
\item $\Ext^n_A(S,A)\ne 0$ for some finite dimensional $A$-module $S$ if
and only if $\Ext^n_A(S,A)\ne 0$ for all finite dimensional $A$-modules
$S$.
\end{enumerate}
\end{enumerate}
Note that (H4) is not left-right symmetric. Hence, when we come to
combine (H4) with duality arguments, we shall need to apply right module
versions of results from Sections \ref{xxsec1}--\ref{xxsec3}.

We next record that Hopf algebras satisfy (H4), along with some other
information. Note that if
$A$ is a Hopf algebra, the base field $k$ can be identified with the
trivial module
$A/(\ker \epsilon)$, where
$\epsilon:  A\to k$ is the counit of $A$.

\begin{lemma}
\label{xxlem4.1}
Let $A$ be a Hopf algebra.
\begin{enumerate}
\item
\cite[Lemma 1.11]{BG1}
For each nonnegative integer $n$, the following are equivalent:
\begin{enumerate}
\item[(1)]
$\Ext^n_A(k,A)\neq 0$.
\item[(2)]
$\Ext^n_A(S,A)\neq 0$ for some finite dimensional
$A$-module $S$.
\item[(3)]
$\Ext^n_A(S,A)\neq 0$ for all nonzero finite dimensional
$A$-modules $S$.
\end{enumerate}
\item
\cite[Lemma 4.8]{WZ1} 
The functors $\Ext^n_A(-,A)$, for $n\ge 0$, are exact on finite
dimensional $A$-modules.
\item 
As an algebra, $A$ satisfies \textup{(H4)}.
\item
If $M$ is a finite dimensional $A$-module, then 
$$\dim_k \Ext^n_A(M,A)= (\dim_k M) \bigl( \dim_k\Ext^n_A(k,A) \bigr)$$ 
for all $n\ge 0$.
\end{enumerate}
\end{lemma}

\begin{proof} Part (c) follows from (a),(b), and (d) follows from the
proof of 
\cite[Lemma 4.8]{WZ1}.
\end{proof}

Recall that an algebra $A$ (respectively, a bimodule $M$) is called {\it
residually  finite dimensional\/} if the intersection of all co-finite
dimensional ideals of $A$ (respectively, sub-bimodules of $M$) is zero.
Such algebras and bimodules appear in our context as follows.

\begin{lemma} \label{newlem4.2}
If $A$ is an algebra satisfying \textup{(H1,H3)}, then all noetherian
$(A,A)$-bimodules are residually finite dimensional. In particular,
$A$ is residually finite dimensional.
\end{lemma}

\begin{proof} Let $M$ be a noetherian
$(A,A)$-bimodule and $N\subseteq M$ a nonzero sub-bi\-mod\-ule. We show
that
$M$ has a co-finite dimensional sub-bimodule not containing $N$.

Set $J := \lann_A(N)$. By \cite[Lemma 8.1]{GW}, $A/J$ embeds 
in $N^{\oplus t}$ as
left $A$-mod\-ules, for some $t\in \NN$. Choose a  maximal ideal
$\mm\supseteq J$, and note that the quotient map $A/J \to A/\mm$ extends
to a nonzero homomorphism $N^{\oplus t} \to E$, where $E$ is the injective hull
of the left $A$-module $A/\mm$. Hence, there exists a nonzero
homomorphism $f: {}_AN \to E$, which extends to a homomorphism $g: {}_AM
\to E$.

Now $A/\mm$ is finite dimensional by (H3), and so $E$ is locally
finite dimensional by Lemma \ref{xxlem3.3}. Hence, $g(M)$ is finite
dimensional, and so the ideal $\qq := \lann_A g(M)$ has finite
codimension in $A$. Since $g(\qq M) =0$, we conclude that $\qq M$ is a
co-finite dimensional sub-bimodule of $M$ that does not contain $N$.
\end{proof} 

\begin{lemma}
\label{xxlem4.2} Suppose $A$ is a residually finite dimensional algebra,
and $M$ a finitely generated, relatively projective $A$-module. Assume
that 
$M/\mm M\neq 0$, or equivalently, $\Hom_A(M,A/\mm)\neq 0$, for 
all co-finite dimensional maximal ideals $\mm$ of $A$. 
Then $M$ is faithful over $A$. If, in addition, $M$ is a 
noetherian $(A,B)$-bimodule for some noetherian ring $B$, then 
$\GKdim M=\GKdim A$.
\end{lemma}

\begin{proof} If, on the contrary, $M$ is not faithful, then 
$\pp:= \lann_A(M)\neq 0$. Since 
$A$ is residually finite dimensional, it has 
a co-finite dimensional ideal $\qq$ such that 
$\pp\not\subseteq \qq$. Consider the short exact sequence
$$0\to (\pp+\qq)/\qq\to A/\qq\to A/(\pp+\qq) \to 0$$
and the associated sequence
$$0\to \bigl( (\pp+\qq)/\qq \bigr) \otimes_A M \to (A/\qq)\otimes_A M
\xrightarrow{f} \bigl( A/(\pp+\qq) \bigr)\otimes_A M\to 0,$$
which is exact because $M$ is relatively projective (Lemma
\ref{xxlem3.4}). The map $f$ is an isomorphism since $\pp M=0$, and so
$\bigl( (\pp+\qq)/\qq \bigr) \otimes_A M =0$. Because
$\pp\not\subseteq \qq$, there is an epimorphism $\bigl( (\pp+\qq)/\qq
\bigr) \to S$ for some simple finite dimensional right $A$-module $S$,
whence $S\otimes_A M =0$. But then $\mm= \rann_A(S)$ is a co-finite
dimensional maximal ideal of $A$ with $(A/\mm) \otimes_A M =0$, which
yields a contradiction to our hypotheses. Therefore $M$ is faithful over
$A$. 

If $M$ is a noetherian bimodule, then $A$ can be embedded
in $M^{\oplus t}$ for some $t$. Hence $\GKdim M=\GKdim A$.
\end{proof}

\begin{lemma} \label{newlem4.4}
Let $A$ be a noetherian algebra with $n := \GKdim A <\infty$, and assume
that $A$ has an Auslander, Cohen-Macaulay dualizing complex $R$.
\begin{enumerate}
\item
$\HB^i(R)=0$ for all $i>0$ and all $i<-n$. 
\item
$\HB^{-n}(R) \ne 0$, and all nonzero left or right $A$-submodules
of $\HB^{-n}(R)$ have GK-dimension $n$.
\item
$\GKdim \HB^{-i}(R) \leq i$ \textup{(}on each side\textup{)} for $0\le
i\le n$.
\end{enumerate}
\end{lemma}

\begin{proof} Note that the Cohen-Macaulay hypothesis implies that $n=
-j_R(A)$, and so $n$ is an integer. Since $A$ is a projective module,
$\RHom_A(A,R) \cong \Hom_A(A,R) \cong R$. Hence,
$$\Ext^q_A(A,R) \cong \HB^q(R)$$
for all $q\in \ZZ$. Since $j_R(A)= -n$, it follows that $\HB^{-n}(R) \ne
0$ and $\HB^i(R)=0$ for all $i<-n$.

If $\HB^i(R) \ne 0$ for some $i$, then because $\HB^i(R)$ is a noetherian
bimodule, the Auslander hypothesis implies that $j_R \bigl( \HB^i(R)
\bigr) \ge i$. Combined with the Cohen-Macaulay condition, this yields
$$\GKdim \HB^i(R) = - j_R \bigl( \HB^i(R) \bigr) \le -i,$$
and similarly on the right.
Part (c) follows, as well as the inequality $0\le -i$. Hence,
$\HB^i(R)=0$ for all $i>0$, and part (a) is proved.

By \cite[Theorem 2.10]{YZ1}, the canonical dimension $\Cdim_R$ is an
exact dimension function. Since
$$\Cdim_R(M)= -j_R(M)= \GKdim M$$
for all finitely generated $A$-modules $M$, it follows that $\GKdim$ is
an exact dimension function. Thus, \cite[Theorem 2.14(1)]{YZ1} implies
that
$\HB^{-n}(R)_A$ is $\GKdim$-pure of GK-di\-men\-sion $n$. By symmetry,
$_A\HB^{-n}(R)$ is $\GKdim$-pure, completing the proof of (b). 
\end{proof}

Recall from Morita theory that if $A$ is a ring, $\Omega$ a progenerator
in $A^\circ\modleft$, and $B :=
\End_{A^\circ}(\Omega)$, then $\Omega$ is an invertible
$(B,A)$-bimodule (with inverse $\Hom_A(\Omega,A)$).

\begin{proposition}
\label{xxprop4.3} Suppose $A$ is an algebra satisfying
\textup{(H1,H2,H3,H4)}, with $n := \GKdim A <\infty$. 
\begin{enumerate}
\item
$A$ has a rigid dualizing complex $R$, and $\HB^i(R) =0$ for all $i\ne
-n$.
\item
$\HB^{-n}(R)$ is an invertible $(A,A)$-bimodule, and $R\cong
\HB^{-n}(R)[n]$.
\item
$\Ext^i_A(S,A)=0$ for all $i\ne n$ and all finite dimensional
$A$-mod\-ules $S$.
\item
If $A$ is a Hopf algebra of finite global dimension, then 
$\gldim A = n$.
\end{enumerate}
\end{proposition}

\begin{proof} (a) By Proposition \ref{xxlem2.1}(a), $A$ has an Auslander,
Cohen-Macaulay, rigid dualizing complex $R$. By 
\cite[Proposition 8.2]{VdB}, rigid dualizing complexes over $A$
are unique up to isomorphism. 

Set $n_0= \max\{ i\in\ZZ \mid \HB^i(R) \ne 0 \}$ and $\Omega=
\HB^{n_0}(R)$. By Lemma \ref{newlem4.4},
$-n\le n_0\le 0$. There is a truncated complex $\Rtil$, quasi-isomorphic to
$R$, such that $\Rtil^i=0$ for all $i>n_0$. For any right $A$-mod\-ule
$M$, we may take an injective resolution $I$ of $M$ and compute
$\Ext^i_{A^\circ}(R,M) \cong \Ext_{A^\circ}(\Rtil,M) \cong \HB^i
\Hom_{A^\circ}(\Rtil,I)$ for all
$i$. Since $I$ vanishes in negative degree, the complex
$\Hom_{A^\circ}(\Rtil,I)$ vanishes in degrees less than $-n_0$, and hence
$\Ext^i_{A^\circ}(R,M) =0$ for all $i< -n_0$. In degrees
$-n_0$ and $-n_0+1$, the complex $\Hom_{A^\circ}(\Rtil,I)$ has the form
$$\Hom_{A^\circ}(\Rtil^{n_0},I^0) \to\Hom_{A^\circ}(\Rtil^{n_0},I^1)
\oplus \Hom_{A^\circ}(\Rtil^{n_0-1},I^0),$$
from which we see that $\Ext^{-n_0}_{A^\circ}(R,M) \cong
\Hom_{A^\circ}(\Omega,M)$.

We now restrict attention to finite dimensional right $A$-modules,
denoted by
$S$. By duality and Lemma \ref{xxlem2.2}(a),
\begin{equation}
\label{E4.1}
\Ext^i_{A^\circ}(R,S) \cong \Ext^i_A(S',D(R)) \cong \Ext^i_A(S',A)
\tag{E4.1}
\end{equation}
for all $i$, where $S'= \Hom_k(S,k)$. Since these isomorphisms are
natural, and since $\Ext^{-n_0}_A(-,A)$ is exact on finite dimensional
$A$-modules by (H4), we see that $\Hom_{A^\circ}(\Omega,-) \cong
\Ext^{-n_0}_{A^\circ}(R,-)$ is exact on finite dimensional right
$A$-modules. Hence,
$\Omega$ is relatively projective as a right $A$-module, and thus
projective, by Proposition \ref{xxprop3.2}.

Because of Lemma \ref{newlem4.2}, $\Omega$ is residually finite
dimensional, and so $\Hom_{A^\circ}(\Omega,S_0)$ is nonzero for some
finite dimensional right $A$-module $S_0$. Since
$$\Hom_{A^\circ}(\Omega,S) \cong \Ext^{-n_0}_{A^\circ}(R,S) \cong
\Ext^{-n_0}_A(S',A)$$ 
for all finite dimensional right $A$-modules $S$, it
follows from (H4) that
$\Hom_{A^\circ}(\Omega,S)$ is nonzero for all finite dimensional right
$A$-modules
$S$. Proposition \ref{xxprop3.2} and Lemmas \ref{newlem4.2},
\ref{xxlem4.2} now imply that
$\Omega_A$ is a progenerator and $\GKdim \Omega_A= n$. By Lemma
\ref{newlem4.4}(c), $n_0= -n$, and thus $\HB^i(R) =0$ for all $i> -n$.
Combined with Lemma \ref{newlem4.4}(a), this establishes part (a).

(b) Since $R$ has nonzero cohomology only in degree $-n$, it is
quasi-isomorphic to $\HB^{-n}(R)[n]= \Omega[n]$.

We have already proved that $\Omega$ is a progenerator for
$A^\circ\modleft$. By definition of a dualizing complex, the canonical
algebra homomorphism
$$A \to \RHom_{A^\circ}(R,R) \cong
\RHom_{A^\circ}(\Omega[n]_A,\Omega[n]_A) \cong
\End_{A^\circ}(\Omega_A)$$ 
is an isomorphism. Therefore, it follows from Morita theory that $\Omega$
is invertible as an $(A,A)$-bimodule.

(c) If $M$ is a finitely generated right $A$-module and $I$ an injective
resolution of $M$, then
$$\RHom_{A^\circ}(R,M) \cong \RHom_{A^\circ}(\Omega[n],M) \cong
\Hom_{A^\circ}(\Omega[n],I).$$ 
The latter complex vanishes in degrees
less than $n$, and is exact in degrees greater than
$n$, because $\Omega_A$ is projective. Hence, $\Ext^i_{A^\circ}(R,M) =0$
for all $i\ne n$. It now follows from (E4.1) that $\Ext^i_A(T,A) =0$ for
all $i\ne n$ and all finite dimensional $A$-modules $T$.

(d) If $A$ is a Hopf algebra with finite global dimension $d$, then the
trivial $A$-mod\-ule $k$ has projective dimension $d$ \cite[Corollary
1.4]{BG1}. Since $A$ is noetherian, $\Ext^d_A(k,A) \ne 0$, and therefore
$d=n$ by part (c).
\end{proof}

\begin{lemma}
\label{xxlem4.4} Let $A$ be an algebra satisfying 
\textup{(H1,H2)}, and $R$ a rigid dualizing complex over $A$.
If $R\cong \Omega[n]$ where $n=\GKdim A$ and $\Omega$ is an
invertible $(A,A)$-bimodule, then $A$ is Auslander-Gorenstein and 
Cohen-Macaulay, and has a quasi-Frobenius classical quotient ring.
Moreover, $\Spec A$ is catenary, and Tauvel's height formula holds.
\end{lemma}

\begin{proof} By Proposition \ref{xxlem2.1}(a), $A$ has an Auslander,
Cohen-Macaulay, rigid dualizing complex, and this complex must be
quasi-isomorphic to $R$ by uniqueness of rigid dualizing complexes
\cite[Proposition 8.2]{VdB}. Since $R \cong \Omega\otimes_A A[n]$,
the ring $A$ is Auslander-Gor\-en\-stein 
and Cohen-Macaulay by \cite[Proposition 4.3]{YZ2}. Catenarity of $\Spec
A$ follows from Theorem \ref{xxlem2.4}. By \cite[Theorem
6.1(2)(3)]{AjSZ}, $A$ has a quasi-Fro\-ben\-i\-us  classical quotient
ring, and
$\GKdim A/\pp=\GKdim A$ for all minimal primes $\pp$ of $A$. Thus, by
Corollary
\ref{xxlem2.5}, Tauvel's height formula holds.
\end{proof}

\begin{theorem}
\label{xxthm4.5} Let $A$ be an algebra satisfying 
\textup{(H1,H2,H3,H4)}. 
\begin{enumerate}
\item
 $A$ is Auslander-Gorenstein and 
Cohen-Macaulay, and has a quasi-Frobenius classical quotient ring.
\item
$\Spec A$ is catenary, and Tauvel's height formula holds.
\item
There exists an invertible $(A,A)$-bimodule $\Omega$ such that
$\Omega[n]$ is a rigid dualizing complex over $A$, where $n= \GKdim A$.
\item
If $A$ is regular, then $\gldim A= \GKdim A$.
\end{enumerate}
\end{theorem}

\begin{proof} Recall from the discussion at the beginning of Section
\ref{xxsec2} that (H2) implies $\GKdim A <\infty$. Parts (a), (b), (c)
follow from Proposition
\ref{xxprop4.3} and Lemma \ref{xxlem4.4}. 

(d) The Auslander-regular and Cohen-Macaulay conditions imply that
$\gldim A \le \GKdim A$. There exists a nonzero finite dimensional
$A$-module $S$ because of (H3), and if $d= \pdim_A S$, then
$\Ext^d_A(S,A) \ne 0$. Proposition \ref{xxprop4.3}(c) implies that $d=
\GKdim A$, and thus $\GKdim A \le \gldim A$.
\end{proof}

\begin{corollary}
\label{newcor4.8}
If $A$ is an affine PI Hopf algebra satisfying \textup{(H2)}, then
statements
\textup{(a)--(d)} of \textup{Theorem \ref{xxthm4.5}} hold.
\end{corollary}

\begin{proof} The affine PI hypothesis implies that $A$ satisfies
(H1,H3), and (H4) follows from Lemma \ref{xxlem4.1}(c). 
\end{proof}

The corollary recovers \cite[Theorem 0.3]{WZ1}.

Theorem \ref{xxthm0.2} is an immediate consequence of 
Theorem \ref{xxthm4.5}, as follows.

\begin{proof}[Proof of Theorem \textup{\ref{xxthm0.2}}]
Let $A$ be a Hopf algebra satisfying (H1,H2,H3). 
By Lemma \ref{xxlem4.1}(c), $A$ satisfies
(H4). Thus, the desired conclusions follow from Theorem \ref{xxthm4.5}.
\end{proof}

\section{Proof of Theorem \ref{xxthm0.1}}

We say an algebra $A$ of finite GK-dimension satisfies
\begin{enumerate}
\item[(H5)]
 provided that, for every maximal ideal $\pp$ of $A$ and every
nonnegative integer $i$, one has 
$$\Ext^i_A(A/\pp,A) \ne 0$$ 
if and only if $i= \GKdim A-\GKdim A/\pp$.
\end{enumerate}
Note that by our conventions, $\Ext^i_A(A/\pp,A)$ refers to an Ext-group
of the \emph{left} $A$-mod\-ules $A/\pp$ and $A$. 

Condition (H5) is
necessary for our main result, as the following lemma shows.

\begin{lemma} \label{newlem5.1}
If $A$ is an Auslander-Gorenstein, Cohen-Macaulay algebra with finite
\textup{GK}-dimension, then
$A$ satisfies \textup{(H5)}.
\end{lemma}

\begin{proof} Let $n := \GKdim A$, let $\pp$ be a maximal ideal of $A$,
and set $d :=
\GKdim A/\pp$. The Cohen-Macaulay condition implies that $j(A/\pp)= n-d$,
and thus we have $\Ext^{n-d}_A(A/\pp,A) \ne 0$ and $\Ext^i_A(A/\pp,A) =0$
for all $i< n-d$.

Suppose that, for some $i$, the $(A/\pp,A)$-bimodule $N :=
\Ext^i_A(A/\pp,A) \ne 0$. Since $N$ is finitely generated as a right
$A$-module,
$j(N_A) \ge i$ by the Auslander condition. Combined with the
Cohen-Mac\-au\-lay condition, we thus find that $\GKdim (N_A) \le n-i$.
By \cite[Proposition 8.3.14(ii)]{MR}, $\GKdim (_{A/\pp}N) \le \GKdim
(N_A)$ (this result does not require $N$ to be finitely generated on the
left). Since $A/\pp$ is a simple ring, $N$ must be faithful as a left
$(A/\pp)$-mod\-ule. The finite generation of $N$ on the right then implies
that
$A/\pp$ embeds in $N^{\oplus t}$, as left $(A/\pp)$-modules, for some
$t\in\NN$. Consequently,
$\GKdim(_{A/\pp}N) =d$, whence $d\le n-i$. Thus $i\le n-d$, and
therefore $i=n-d$ by the previous paragraph.
\end{proof}

\begin{theorem}
\label{xxthm5.1} Let $A$ be an algebra satisfying \textup{(H1,H2,H5)}.
\begin{enumerate}
\item There exists an invertible $(A,A)$-bimodule $\Omega$ such that
$\Omega[n]$ is a rigid dualizing complex over $A$, where $n= \GKdim A$.
\item $A$ is Auslander-Gorenstein and 
Cohen-Macaulay, and has a quasi-Frobenius classical quotient ring.
\item $\Spec A$ is catenary, and Tauvel's height formula holds.
\end{enumerate}
\end{theorem}

\begin{proof} (a) By Proposition \ref{xxlem2.1}(a), $A$ has an Auslander,
Cohen-Macaulay, rigid dualizing complex $R$.

Fix a maximal ideal $\pp$ of $A$. Let 
$m :=\GKdim A/\pp$ and $n :=\GKdim A$, and let $M_\pp$ denote
the complex over $A/\pp$ given by
$\RHom_A(A/\pp,R)$. By Proposition \ref{xxlem2.1}(b),
$M_\pp$ is an Auslander, Cohen-Macaulay, rigid dualizing complex over
$A/\pp$.  Since $A/\pp$ is simple, it follows from \cite[Theorem
0.2]{YZ2} that
$M_\pp=\Omega_\pp [m]$ where
$\Omega_\pp$ is an invertible $(A/\pp,A/\pp)$-bimodule. By duality,
$$\Ext^i_{A^\circ}(R[-n],\Omega_\pp)=
\Ext^{i+n-m}_{A^\circ}(R,M_\pp)\cong 
\Ext^{i+n-m}_A(A/\pp,A),$$
which is nonzero if and only if $i=0$ (by hypothesis (H5)). 
By $\Hom$-$\otimes$ adjunction in $\D^b({A^\circ}\modleft)$, we
have
$$\Ext^i_{A^\circ}(R[-n],\Omega_\pp)\cong 
\Hom_{(A/\pp)^\circ}(\Tor_i^A(R[-n],A/\pp),\Omega_\pp).$$
Since $\Tor^A_i(R[-n],A/\pp)$
is projective over $A/\pp$ (Lemma \ref{xxlem1.6}(b)), 
it follows that $\Tor_i^A(R[-n],A/\pp)\neq 0$ if and only if
$i=0$. 

By Proposition \ref{xxprop1.1}, $R[-n]$ is isomorphic to a
bimodule, say $\Omega$, and $\Omega$ is a progenerator
for ${A^\circ}\modleft$. As in the proof of Proposition
\ref{xxprop4.3}(b), the natural map $A\to \End_{A^\circ}(\Omega)$ is an
isomorphism, and Morita theory implies that
$\Omega$ is invertible.

(b)(c) These properties follow from part (a) and Lemma 
\ref{xxlem4.4}.
\end{proof}

For future use, we record the following consequence of Lemma
\ref{newlem5.1} and Theorem \ref{xxthm5.1}.

\begin{corollary} \label{newcor5.3}
If $A$ is an Auslander-Gorenstein, Cohen-Macaulay algebra satisfying
\textup{(H1,H2)}, then there exists an invertible $(A,A)$-bimodule
$\Omega$ such that
$\Omega[n]$ is a rigid dualizing complex over $A$, where $n= \GKdim A$.
\end{corollary}

In the rest of this section, we verify the 
hypothesis (H5) for standard generic quantized coordinate rings $\OqG$.
Since this relies on results from \cite{Jos}, we will need to assume
that $q$ is transcendental over $\QQ$.

\begin{lemma} 
\label{lemma1} Suppose $A_1,\dots,A_n$ are  connected
graded noetherian algebras having enough  normal elements.
\begin{enumerate}
\item The algebra $A := A_1\otimes\cdots\otimes A_n$ is
connected graded noetherian and has enough normal elements.
\item If $\sigma$ is a graded automorphism of $A$, then 
$A[t;\sigma]$  is connected graded noetherian with enough normal
elements, where $t$ is assigned degree $1$.
\item Let $C:= A[t_1;\sigma_1]\cdots[t_m;\sigma_m]$ be an
iterated skew polynomial ring, where the $\sigma_i$ are graded
automorphisms and the $t_i$ have degree $1$. Then $C$ is
connected graded, universally  noetherian, and has enough normal
elements.
\item If the $A_i$ are Gorenstein \textup{(}respectively,
regular\textup{)},  then the algebra $C$ above is
Auslander-Gorenstein \textup{(}respectively,
Auslander-regular\textup{)} and Cohen-Mac\-aulay.
\end{enumerate}
\end{lemma}

\begin{proof} (a) It follows from 
\cite[Propositions 4.3, 4.9 and 4.10]{ArSZ} that  every
connected graded noetherian algebra with enough normal elements
is universally noetherian. As a consequence, $A$ is universally
noetherian by induction on $n$. It is clear that $A$ is
connected graded (with respect to the standard tensor product
grading).  To show that $A$ has enough normal elements, we can
proceed by induction on $n$, so it suffices to consider the case
$n=2$. Let $P$ be a non-maximal graded prime ideal of
$A_1\otimes A_2$. Since $A_1\otimes 1$ commutes with $1\otimes
A_2$,  the preimage $Q_i$ of $P$ in $A_i$ is prime. Since we can
pass to the algebra $A/( Q_1\otimes A_2+ A_1\otimes Q_2)$, there
is no loss of generality in assuming that $A_1$ and $A_2$ are
graded prime, that $P\cap (A_1\otimes 1) =0$, and that $P\cap
(1\otimes A_2) =0$. Since $A/P$ is not simple, $A$ is not finite
dimensional, so $A_1$ and $A_2$ cannot both be finite
dimensional. Say $A_1$ is infinite dimensional, so that it is
not simple. Then there is a homogeneous normal element $x_1\in
A_1$ with positive degree, whence the image of $x_1\otimes 1$ in
$A/P$ is a homogeneous normal element, with positive degree
because it is nonzero. The proof is symmetric in case $A_2$ is
infinite dimensional. This shows that $A_1\otimes A_2$ has
enough normal elements.  

(b) It is clear that $A[t;\sigma]$ is connected graded noetherian. Let
$P$ be a non-max\-i\-mal graded prime ideal of $A[t;\sigma]$. If
$t\notin P$, then the image of
$t$ in $A[t;\sigma]/P$ is a homogeneous normal element of degree 1.
Otherwise, $A[t;\sigma]/P \cong A/P'$ for some non-maximal graded prime
ideal $P'$ of $A$. The assertion follows.  

(c) Using (a), (b), and induction, we see that $C$ is connected graded
noetherian with enough normal elements. The universally noetherian
property follows from the comment at the beginning of the proof of (a).

(d) Recall that the global dimension of a connected graded noetherian
algebra equals the projective dimension of its trivial module (e.g.,
\cite[Chapter 1, Corollary 8.7]{NVO}). For $A_1$ and $A_2$, this means 
\begin{equation}
\begin{aligned}
\gldim A_1\otimes A_2 &=\pdim_{A_1\otimes A_2} k \\
 &\leq  \pdim_{A_1}k +\pdim_{A_2}k=
\gldim A_1+\gldim A_2.
\end{aligned}  \notag
\end{equation}
Hence, if $A_1$ and $A_2$ are regular, then so is $A_1\otimes A_2$.

Next, we assume that $A_1$ and $A_2$ are Gorenstein
of injective dimensions $d_1$ and $d_2$ respectively. 
Since $A_1$ and $A_2$ have enough normal 
elements, they are Artin-Schelter 
Gorenstein \cite[Theorems 0.2 and 0.3]{Zh} 
in the following sense:
$$\Ext^i_{A_j}(k,A_j)=
\begin{cases} k & {\rm{if}}\;\; i=d_j\\ 0 & 
{\rm{if}}\;\; i\neq d_j
\end{cases}$$
for $j=1,2$.
It follows from the K{\"u}nneth formula that 
\begin{equation}
\begin{aligned}
\Ext^{n}_{A_1\otimes A_2}(k,A_1\otimes A_2)
 &=\bigoplus_{i=0}^n \Ext^i_{A_1}(k,A_1)\otimes
\Ext^{n-i}_{A_2}(k,A_2)  \\
 &=\begin{cases} k& {\rm{if}}\;\; n=d_1+d_2\\ 0& 
{\rm{if}}\;\; n\neq d_1+d_2
\end{cases}.
\end{aligned}  \notag
\end{equation}
This means that ${\rm{k.}}\injdim 
A_1\otimes A_2=d_1+d_2$, where ${\rm{k.}}\injdim$
is the $k$-injective dimension defined in
\cite{Jor}. By \cite[Corollary 8.12]{AZ1}, $A_1\otimes A_2$
satisfies the condition $\chi$, and its cohomological dimension 
$\cd(A_1\otimes A_2)$ is bounded by its Krull dimension. Since
$A_1\otimes A_2$ has enough normal elements, its Krull dimension is
finite \cite[Proposition 0.9]{YZ1}, and thus $\cd(A_1\otimes A_2) <
\infty$. Hence, we can apply 
\cite[Theorem 4.5]{Jor} to obtain that 
$\injdim A_1\otimes A_2={\rm{k.}}\injdim 
A_1\otimes A_2$. Therefore $A_1\otimes A_2$
is Gorenstein.

It follows by induction that 
if the $A_i$ are Gorenstein (respectively, regular), then so
is $A$. Either property then passes to $C$ (e.g., use
\cite[Theorem 3.6(1)]{Lev} and \cite[Theorem 5.3(iii)]{MR}
inductively). The assertion now follows from (c) and
\cite[Theorem 0.2]{Zh}.
\end{proof}

\begin{lemma}
\label{lemma2} Let $A$ be a noetherian connected graded 
Gorenstein \textup{(}respectively, regular\textup{)} algebra
with  enough normal elements.  Let $C := A [t_1;\sigma_1]\cdots 
[t_m;\sigma_m]$ be an iterated skew polynomial ring, where the
$\sigma_i$ are graded  automorphisms and the $t_i$ have degree
$1$. Let $D := C[t_1^{-1},\dots,t_m^{-1}]= A
[t_1^{\pm1};\sigma_1]\cdots  [t_m^{\pm1};\sigma_m]$. Suppose
that $t_j$ is an eigenvector for $\sigma_i$, for all $i>j$. Then:
\begin{enumerate}
\item
$D$ is Auslander-Gorenstein \textup{(}respectively, 
Auslander-regular\textup{)} and Cohen-Macaulay.
\item
$D$ has a filtration such that $\gr D$ is a noetherian,
connected graded algebra with enough normal elements.
\end{enumerate}
\end{lemma}

\begin{proof} Assume that $A$ is regular. The Gorenstein case is
completely parallel, except that one has to check the Gorenstein
case of \cite[Lemma]{LS}, to feed into \cite[Lemmas I.15.4,
II.9.11]{BG2}.

(a) By Lemma \ref{lemma1}, $C$ is connected graded noetherian,
Auslander-regular, and Cohen-Macaulay. Observe that each $t_i$
is regular and normal in $C$; moreover, $t_iC_j= C_jt_i$ for all
$i,j$, where the $C_j$ are the homogeneous components of $C$.
Hence, the element $t := t_1t_2\cdots t_m$ is regular normal,
and $tC_j= C_jt$ for all $j$. By \cite[Lemmas II.9.11(b)]{BG2},
$C[t^{-1}]= D$ is Auslander-regular and Cohen-Macaulay. 

(b) Let $t$ be as above; then there exists a graded automorphism
$\sigma$ of $C$ such that $tc= \sigma(c)t$ for all $c\in C$. Set
$E := C[T;\sigma^{-1}]$, which is a connected graded noetherian
algebra with enough normal elements, by Lemma \ref{lemma1},
where $\deg(T)=1$. Give $E$ the filtration induced from its
grading. Since $E/ \langle tT-1\rangle \cong D$, the required
properties pass from $E$ to $D$.
\end{proof}

\begin{theorem} \label{newthm5.4}
Let $G$ be any connected semisimple Lie group
over $\CC$, and $A := \OqG$ the standard quantized coordinate
ring. Assume that $q$ is transcendental over $\QQ$. There exists
a left and right Ore set $\calC \subset A$ such that
\begin{enumerate}
\item $(A/I)\calC^{-1} \ne 0$ for all proper ideals $I$ of $A$.
\item The algebra $B := A\calC^{-1}$ is Auslander-regular and
Cohen-Macaulay.
\item $B$ satisfies \textup{(H2)}.
\item $\GKdim B = \dim G$.
\end{enumerate}
\end{theorem}

\begin{proof} We first work over the field $F := \QQ(q) \subset
\CC$ and then extend scalars. Let $\fg$ be the Lie algebra of
$G$, and write $\Ucheck_F$ for the simply connected quantized
enveloping algebra of $\fg$ over $F$. Following the notation of
\cite[\S I.6.3]{BG2}, we write $E_i$, $F_i$, $K_\lambda^{\pm1}$
for the standard generators of $\Ucheck_F$, where $i= 1,\dots,
n= \rank(\fg)$ and $\lambda$ runs through the weight lattice
$P$. Write $U_F^+$ and $U_F^-$ for the respective subalgebras of
$\Ucheck_F$ generated by the $E_i$ and the $F_i$. Since the
quantum Serre relations for the $E_i$ and the $F_i$ are
homogeneous, $U_F^+$ and $U_F^-$ are connected graded
$F$-algebras, where the $E_i$ and $F_i$ are homogeneous elements
of degree 1.

Write $A_F$ for the version of $\OqG$ defined over $F$, and let
$\calC$ denote the multiplicative subset of $A_F$ generated by
the coordinate functions $c^{\mu}_{\mu,\mu}$ (in the notation of
\cite[\S9.1.1]{Jos}), for $\mu\in P^+$. By \cite[Lemma
9.1.10]{Jos}, $\calC$ is an Ore set in $A_F$. Moreover,
$$A_F\calC^{-1} \cong (U_F^+\otimes_F
U_F^-)[K_\lambda^{-1}\otimes K_\lambda \mid \lambda\in P]
\subset \Ucheck_F\otimes_F \Ucheck_F$$
\cite[Lemma 9.2.13, Proposition 9.2.14]{Jos}, and $\calC$ is
disjoint from all prime ideals of $A_F$ \cite[Corollary
9.3.9]{Jos}. It follows that $\calC$ is disjoint from all proper
ideals $I$ of $A_F$, so that $(A_F/I)\calC^{-1} \ne 0$.

We now identify $A$ with $A_F\otimes_F \CC$ and $A_F$ with the
$F$-subalgebra $A_F\otimes 1$. Then $\calC$ becomes an Ore set
in $A$, and since every proper ideal of $A$ contracts to a
proper ideal of $A_F$, we see that (a) holds. Extending scalars,
the isomorphism above yields
$$B= A\calC^{-1} \cong B' := (U^+\otimes_{\CC} U^-)
[K_\lambda^{-1}\otimes K_\lambda \mid \lambda\in P] \subset
\Ucheck\otimes_{\CC} \Ucheck,$$ where $U^{\pm} := U^{\pm}_F
\otimes_F \CC$ and $\Ucheck := \Ucheck_F \otimes_F \CC$.

As with $U^{\pm}_F$, the $\CC$-algebras $U^{\pm}$ are connected
graded and noetherian. The basic relations for $\Ucheck$ show
that the inner automorphism induced by any $K_\lambda$ on
$\Ucheck$ restricts to graded automorphisms of $U^+$ and $U^-$.
Since $P$ is a free abelian group of rank $n$, we see that $B'$
is an iterated skew-Laurent extension of the form
\begin{equation}
B'= (U^+\otimes_{\CC} U^-) [t_1^{\pm1};\sigma_1] \cdots
[t_n^{\pm1};\sigma_n],
\tag{E5.1}
\end{equation} 
for some graded automorphisms
$\sigma_i$, where $\deg(t_j) =1$ for all $j$ and $\sigma_i(t_j)=
t_j$ for $i>j$. Therefore, properties (b) and (c) will follow
from Lemmas \ref{lemma1} and \ref{lemma2} once we know that
$U^+$ and $U^-$ are regular and have enough normal elements.
Since $U^+\cong U^-$, we need only deal with $U^+$.

Ringel has shown that $U^+_F$ is isomorphic to the twisted
generic Hall algebra $\cal H$ \cite{Rin1} (see also
\cite[Theorem 3]{Gr}), and that
$\cal H$ is an iterated skew polynomial ring of the form
$${\cal H}= F[X_1][X_2;\tau_2,\delta_2] \cdots
[X_m;\tau_m,\delta_m]$$  
where the following hold \cite[Theorems 2,3]{Rin2}:
\begin{enumerate}
\item[(1)] The number $m$ of iterations equals the number of positive
roots for the Lie algebra $\fg$ of $G$.
\item[(2)] Each $\tau_i$ and each $\delta_i$ is
$F$-linear.
\item[(3)] $X_j$ is an eigenvector of $\tau_i$ for all $i>j$.
\item[(4)] Each $\delta_i\tau_i= q^{r_i}\tau_i\delta_i$ for some
positive integer $r_i$.
\end{enumerate}
Consequently,
$$U^+ \cong U_F^+ \otimes_F \CC \cong
\CC[X_1][X_2;\tauhat_2,\deltahat_2] \cdots
[X_m;\tauhat_m,\deltahat_m]$$ with corresponding properties.
Hence, $U^+$ is regular \cite[Theorem 7.5.3]{MR}, and all its
prime ideals are completely prime \cite[Theorem 2.3]{GLet}. 

The existence of a large supply of normal elements was
established by Caldero \cite{Cal} for the $\CC(\qhat)$-form of
$U^+$, where $\qhat$ is an indeterminate over $\CC$. Rather than
reworking all his proofs for $U^+$, we can obtain what we need
from his results with some base change arguments.

First, extend $\{q\}$ to a transcendence base $X$ for $\CC$ over
$\QQ$, and let $K$ denote the algebraic closure of
$\QQ(X\setminus \{q\})$ within $\CC$. Then $q$ is transcendental
over $K$, and $K\cong \CC$ (as abstract fields) because $X$ is
infinite. Hence, there is an isomorphism $K(q) \rightarrow
\CC(\qhat)$ sending $q \mapsto \qhat$. Consequently, the algebra
$U^+_{K(q)} := U^+_F \otimes_F K(q)$ is isomorphic, as a ring,
to the $\CC(\qhat)$-form of $U^+$ (that is, the
$\CC(\qhat)$-algebra given by generators $E_1,\dots,E_n$
satisfying the quantum Serre relations). Caldero's result,
\cite[Corollaire 3.2]{Cal}, thus shows that every ideal of
$U^+_{K(q)}$ has a normalizing sequence of generators.

Let $\pp$ be a non-maximal graded prime ideal of $U^+$. Since
$\pp$ is completely prime, the contraction $\pp_1 := \pp\cap
U^+_{K(q)}$ is a graded prime ideal of $U^+_{K(q)}$. Let $\fm_1$
denote the (unique) graded maximal ideal of $U^+_{K(q)}$, and
observe that $\fm_1 U^+$ is the graded maximal ideal of $U^+$.
Consequently, $\fm_1 \nsubseteq \pp$, and so $\pp_1$ is properly
contained in $\fm_1$. Since $\fm_1$ has a normalizing sequence
of generators, there exists a nonzero element $c\in \fm_1/\pp_1$
which is normal in $U^+_{K(q)}/\pp_1$. A leading term argument,
using the fact that $U^+_{K(q)}/\pp_1$ is a domain, shows that
the highest degree term of $c$ is normal, so there is no loss of
generality in assuming that $c$ is homogeneous. Since $c\in
\fm_1/\pp_1$, its degree is positive. Thus, finally, the image
of $c$ in $U^+/\pp$ is a homogeneous normal element of positive
degree.

Therefore $U^+$ has enough normal elements, as required. Statements (b)
and (c) are now proved.

Since $U^-\cong U^+$, the algebra $U^+\otimes_{\CC} U^-$ can be
written as an iterated skew polynomial extension of $U^+$, and so
$$U^+\otimes_{\CC} U^- \cong \CC[X_1][X_2;\tauhat_2,\deltahat_2]
\cdots [X_m;\tauhat_m,\deltahat_m] [Y_1][Y_2;\tautil_2,\deltil_2]
\cdots [Y_m;\tautil_m,\deltil_m],$$
where 
\begin{enumerate}
\item[(5)] Each $\tauhat_i$, $\tautil_i$, $\deltahat_i$, and
$\deltil_i$ is $\CC$-linear.
\item[(6)] $X_j$ is an eigenvector of $\tauhat_i$ for all $i>j$.
\item[(7)] $Y_j$ is an eigenvector of $\tautil_i$ for all $i>j$,
and $X_j$ is fixed by $\tautil_i$ for all $i,j$.
\end{enumerate}
In view of (5)--(7), it follows from \cite[Lemma II.9.7]{BG2}
that
$$\GKdim (U^+\otimes_{\CC} U^-) = 2m.$$
Now return to equation (E5.1), and consider the intermediate iterations
$$B'_j := (U^+\otimes_{\CC}
U^-) [t_1^{\pm1}; \sigma_1] \cdots [t_j^{\pm1}; \sigma_j].$$
Recall that $U^+\otimes_{\CC} U^-$ is generated by the elements
$E_i\otimes 1$ and $1\otimes F_i$, each of which is an eigenvector for the
automorphisms
$\sigma_j$, and that $t_l$ is fixed by $\sigma_j$ when $l<j$. Hence,
$B'_{j-1}$ is generated by a finite set of $\sigma_j$-eigenvectors, from
which we see that $\sigma_j$ is {\it locally algebraic\/} in the sense of
\cite{LMO}, meaning that for each $b\in B'_{j-1}$, the forward orbit $\{
\sigma_j^t(b) \mid t\ge 0\}$ spans a finite dimensional subspace of
$B'_{j-1}$. Consequently, it follows from \cite[Proposition 1]{LMO} that
$\GKdim B'_j = 1+ \GKdim B_{j-1}$.
Hence, 
$$\GKdim B= \GKdim B'= 2m+n= \dim \fg = \dim G$$
(recall (1) above), and part (d) is proved.
\end{proof}

\begin{lemma} \label{newlem5.5}
Let $A$, $\calC$, and $B$ be as in \textup{Theorem \ref{newthm5.4}}.
\begin{enumerate}
\item
Any noetherian $(A,A)$-bimodule $M$ is $\calC$-torsionfree as a left and a
right $A$-module. Hence, if $M\ne 0$, then $M\calC^{-1} \ne 0$ and
$\calC^{-1}M \ne 0$.
\item
$\GKdim (A/\pp)\calC^{-1}= \GKdim A/\pp$ for all prime ideals $\pp$ of $A$.
\item
$\gldim A= \GKdim A= \gldim B= \GKdim B <\infty$.
\item
 Let $n= \gldim A$. If $S$ is any nonzero finite dimensional left or
right $A$-mod\-ule, then
$\Ext^n_A(S,A) \ne 0$, and $\Ext^i_A(S,A) =0$ for all $i\ne n$.
\end{enumerate}
\end{lemma}

\begin{proof} (a) Let $N$ be the left $\calC$-torsion submodule of $M$.
Then $N$ is a sub-bimodule, and so $N= x_1A +\dots+ x_tA$ for some
elements $x_1,\dots,x_t$. There exists $c\in \calC$ such that $cx_i=0$ for
all $i$, and thus $c\in \lann_A(N)$. Since $\calC$ is disjoint from all
proper ideals of $A$ (Theorem \ref{newthm5.4}(a)), it follows that
$\lann_A(N) =A$ and $N=0$. Thus $M$ is $\calC$-torsionfree as a left
$A$-module. By symmetry, it is also $\calC$-torsionfree as a right
$A$-module.

(b) If $\pp$ is a prime ideal of $A$, then by part (a), all elements of
$\calC$ are regular modulo $\pp$. (This also follows directly from the fact
that $\calC\cap \pp = \varnothing$, by \cite[Lemma 10.19]{GW}.) Hence, $\pp
B$ is a prime ideal of $B$ which contracts to $\pp$ \cite[Theorem
10.20]{GW}. Consequently, $A/\pp$ embeds in $B/\pp B \cong
(A/\pp)\calC^{-1}$, and these prime rings have isomorphic Goldie quotient
rings, call them $Q_1= \Fract A/\pp$ and $Q_2= \Fract (A/\pp)\calC^{-1}$.
By Theorems \ref{xxprop2.6} and
\ref{newthm5.4}, $A$ and $B$ satisfy (H2), hence so do $A/\pp$ and
$(A/\pp)\calC^{-1}$. It now follows from \cite[Theorem 6.9(a)]{YZ6} that
$\Htr Q_1= \GKdim A/\pp$ and $\Htr Q_2= \GKdim (A/\pp)\calC^{-1}$. Since
$Q_1\cong Q_2$, part (b) is proved.

(c) We have $\GKdim A= \GKdim B$ by part (b), $\gldim B\le \gldim A$
because $B$ is a localization of $A$, and $\gldim B \le \GKdim B$
because $B$ is Auslander-regular and Cohen-Macaulay.

Set $n := \GKdim B$, and let $k$ denote the trivial module for $A$. We
have $k\calC^{-1} \ne 0$ by part (a), from which we see that $k\calC^{-1} =k$,
that is, the $A$-module structure on $k$ extends to a $B$-module
structure. Since $B$ is Cohen-Macaulay, $j(_Bk) = n$, and thus
$\Ext^n_B(k,B) \ne 0$. Hence, $n\le \pdim_B k\le \gldim B$, and we have
$$\gldim B= \pdim_B k= n.$$
As noted in the proof of \cite[Proposition 2.7]{BG1}, it follows from
the flatness of $B$ over $A$ that $\pdim_A k \le n$, and thus $\gldim
A\le n$ by \cite[Corollary 1.4(c)]{BG1}. Therefore
$$\gldim A= \GKdim A= \pdim_A k= n.$$

(d) By symmetry, it is enough to prove this for left modules. In view of
\cite[Lemma 1.11]{BG1} (i.e., Lemma \ref{xxlem4.1}(a)), it suffices to
consider the case when $S$ is the trivial module $k$. 

We proved above that $\pdim_A k= n$, whence $\Ext^n_A(k,A) \ne 0$, and
also that $j(_Bk) =n$, whence $\Ext^i_B(k,B) =0$ for all $i\ne n$. It
now follows using \cite[Proposition 1.6]{BL} that
$$\Ext^i_A(k,A)\calC^{-1} \cong
\Ext^i_{A\calC^{-1}}(k\calC^{-1},A\calC^{-1}) =0$$ 
for all $i\ne n$. However,
$\Ext^i_A(k,A)$ is finite dimensional by Lemma \ref{xxlem2.3}, and so is
noetherian as an $(A,A)$-bimodule. Thus, it follows from part (a) that 
$\Ext^i_A(k,A) =0$ for all $i\ne n$.
\end{proof}

\begin{lemma} \label{xxlem5.4} 
Let $G$ be any connected semisimple Lie group
over $\CC$, and $A := \OqG$ the standard quantized coordinate
ring. If $q$ is transcendental over $\QQ$, then $A$ satisfies
\textup{(H5)}, and $\GKdim A= \dim G$.
\end{lemma}

\begin{proof} Let $\calC$ and $B$ be as in Theorem \ref{newthm5.4}, and
set
$$n= \gldim A= \GKdim A= \gldim B= \GKdim B= \dim G$$
(recall Lemma \ref{newlem5.5}(c) and Theorem \ref{newthm5.4}(d)). Fix a
maximal ideal
$\pp
\subset A$, and set $d := \GKdim A/\pp$ and $\qq := \pp B$. Then $\qq$ is
a maximal ideal of $B$ by localization theory \cite[Theorem 10.20]{GW},
and
$\GKdim B/\qq= d$ by Lemma \ref{newlem5.5}(b). Since $B$ is
Auslander-regular and Cohen-Macaulay, Lemma \ref{newlem5.1} implies that
for nonnegative integers $i$, we have $\Ext^i_B(B/\qq,B) \ne 0$ if and
only if $i= n-d$.

According to Lemma \ref{xxlem2.2}(b), each $\Ext^i_A(A/\pp,A)$ is a
noetherian
$(A,A)$-bimodule. Moreover, $\Ext^i_A(A/\pp,A) \calC^{-1} \cong
\Ext^i_B(B/\qq,B)$ by \cite[Proposition 1.6]{BL}. Taking account of
Lemma \ref{newlem5.5}(a), we conclude that $\Ext^i_A(A/\pp,A) \ne 0$ if
and only if $i=n-d$, as desired. 
\end{proof}

Now we are ready to prove Theorem \ref{xxthm0.1}.

\begin{proof}[Proof of Theorem \textup{\ref{xxthm0.1}}]
By Theorems \ref{xxlem1.9}, \ref{xxprop2.6} and Lemmas \ref{newlem5.5},
\ref{xxlem5.4},
$\OqG$ satisfies the hypotheses (H1,H2,H5), and $\gldim \OqG= \GKdim
\OqG= \dim G$. Theorem \ref{xxthm0.1}
thus follows from Theorems \ref{xxprop2.6} and
\ref{xxthm5.1}.
\end{proof}

\providecommand{\bysame}{\leavevmode\hbox to3em{\hrulefill}\thinspace}
\providecommand{\MR}{\relax\ifhmode\unskip\space\fi MR }
\providecommand{\MRhref}[2]{%
   \href{http://www.ams.org/mathscinet-getitem?mr=#1}{#2} }
\providecommand{\href}[2]{#2}

\end{document}